\documentclass[12pt]{article} 
\usepackage{latexsym}

\makeatletter

\@addtoreset{equation}{section}
\makeatother

\newtheorem{theorem}{Theorem}[section]
\newtheorem{prop}[theorem]{Proposition}
\newtheorem{lemma}[theorem]{Lemma}
\newtheorem{cor}[theorem]{Corollary}
\newtheorem{remark}[theorem]{Remark}

\newtheorem{claim}[theorem]{Claim}
\newtheorem{sublemma}[theorem]{Sublemma}

\newcommand{\R}{{\bf R}}
\newcommand{\C}{{\bf C}}

\newcommand{\Z}{{\bf Z}}
\newcommand{\qu}{{\bf H}}
\newcommand{\A}{{\cal A}}
\newcommand{\G}{{\cal G}}

\newcommand{\W}{{\cal W}}

\newcommand{\M}{{\cal M}}

\newcommand{\PP}{{\cal P}}
\newcommand{\CC}{{\cal C}}

\newcommand{\V}{{\cal V}}

\newcommand{\ZZ}{{\cal Z}}
\newcommand{\sw}{{\rm sw}}

\newcommand{\XX}{{\cal X}}

\newcommand{\bH}{\mathbf{H}}

\newcommand{\stab}{\mbox{\rm stab}}
\newcommand{\kernel}{{\rm ker\,}}
\newcommand{\coker}{{\rm coker\,}}
\newcommand{\ran}{{\rm Ran}}

\newcommand{\proof}{\textbf{Proof\ }}
\newcommand{\proofof}[1]{\textbf{Proof of #1\ }}

\newcommand{\qed}{ \mbox{qed.}}

\newcommand{\normi}[2]{\|#1\|_{#2}}

\newcommand{\const}{\mbox{\rm const.}}

\newcommand{\dimn}{\mbox{\rm dim}}

\newcommand{\cs}{{\rm cs}}
\newcommand{\csd}{{\rm csd}}
\newcommand{\su}{{\rm su}}
\newcommand{\ad}{{\rm ad}}
\newcommand{\Ad}{{\rm Ad}}

\newcommand{\SF}{{\rm SF}}

\newcommand{\detind}{{\rm detind\,}}
\newcommand{\sym}{{\rm Sym}}

\renewcommand{\epsilon}{\varepsilon}
\newcommand{\mPhi}{{\mit\Phi}}
\newcommand{\mPsi}{{\mit\Psi}}
\newcommand{\mLambda}{{\mit\Lambda}}
\newcommand{\mXi}{{\mit\Xi}}

\newcommand{\mTheta}{{\mit\Theta}}

\begin{document}

\title{Defining an SU(3)-Casson/U(2)-Seiberg-Witten
integer invariant for integral homology 3-spheres}
\author{Yuhan Lim}
\date{}

\maketitle
\begin{abstract}
    \noindent The SU(3)-Casson invariant for integral homology 3-spheres
    as studied by Boden-Herald possesses a \lq spectral flow 
    obstruction\rq\ to being
    an integer valued invariant which 
    depends only on the non-degenerate (perturbed) moduli 
    space of flat SU(3)-connections.
    This obstruction is the non-trivial spectral flow of a family of 
    twisted signature operators in 3-dimensions.
    The parallel U(2)-Seiberg-Witten construction also has an
    obstruction but from the non-trivial spectral 
    flow of a family of twisted Dirac operators.  
    By taking the SU(3)-flat and U(2)-Seiberg-Witten equations simultaneously
    the obstructions can be made to cancel and an integer 
    invariant is obtained. 
\end{abstract}

\setlength{\parskip}{0pt}
\setlength{\parskip}{6pt plus3pt minus 3pt}
\goodbreak

\section{Introduction}

In 1985 Casson \cite{AM} introduced his now well-known integer invariant for 
(oriented) integral homology 3-spheres (ZHS). This beautiful invariant is a lift
of the $\bmod 2$ $\mu$-invariant and with it Casson showed 
how to prove a number of remarkable theorems in low-dimensional 
topology.  Roughly speaking Casson's invariant algebraically counts, up to conjugacy 
the number of representations of the fundamental group into $SU(2)$. 
Shortly after Taubes \cite{taubes} showed how to interpret this as an 
Euler characteristic/Hopf 
index in the infinite dimensional setting of gauge theory, the bridge being the 
correspondence between flat connections on an $SU(2)$-bundle and 
representations into $SU(2)$ of the fundamental group. Meanwhile 
Floer \cite{floer}
defined his homology groups based on Witten's Morse theory ideas 
applied to gauge theory -- the Casson invariant appears as (half) the 
Euler characteristic of the Floer groups.

A natural idea is to extend Casson's invariant by utilizing 
gauge groups different from $SU(2)$, especially to higher $SU(n)$. A 
proposal along Casson's original approach involving Heegaard splittings and 
representation varieties was announced by 
Cappell-Lee-Miller  \cite{CLM1} in 1990. 
Then in 1998 Boden-Herald \cite{BH} presented the 
first detailed account of an $SU(3)$-Casson invariant based on Taubes'
interpretation of Casson's invariant. 

\subsection{Spectral flow obstruction in the $\mathbf{SU(3)}$-Casson invariant}

Let $Y$ denote an oriented ZHS and $F\to Y$ an $SU(3)$-bundle. In the 
$SU(3)$-Casson invariant consider $\A_{F}$ the space of 
$SU(3)$-connections on $F$ and $\G_{F}$ the group of gauge 
transformations of
$F$. (For the purpose of this introduction we omit the Sobolev 
completions of these $C^{\infty}$ objects. Details are in the main 
body of the article.)
The flat connections on $F$ are the critical points of the 
Chern-Simons functional $\cs$ on $\A_{F}$. The quotient space
$\A_{F}/\G_{F}$ is stratified; the highest comes from the \textit{irreducibles} 
$\A^{*}_{F}$ (finite stabilizer under $\G_{F}$) and the lowest the
\textit{trivial} connections $\A^{0}_{F}$. The only intermediate strata relevant 
for us comes from the \textit{$U(2)$-reducibles} $\A^{I}_{F}$ (those 
with $U(1)$-stabilizer).
For a ZHS the moduli space of flat connections splits as the stratas 
above:
\begin{equation}
    \M^{\su}=\M^{\su*}\cup\M^{\su,I}\cup\{[\mTheta]\}
\end{equation}
where $[\mTheta]$ is the orbit of any trivial connection $\mTheta$.
$\M^{\su,I}$ is exactly the moduli space of flat $U(2)$-connections.

A central ingredient is the necessity of perturbing the 
Chern-Simons function in order to make the critical 
points $\M^{\su}$ a finite set of non-degenerate points. It is a 
theorem that this is always possible. Let us 
denote the perturbed moduli space as $\M^{\su}_{\pi}$. Then in the same 
manner as the splitting above we have
\begin{equation}
    \M^{\su}_{\pi}=\M^{\su *}_{\pi}\cup\M^{\su,I}_{\pi}\cup\{[\mTheta]\}.
\end{equation}
Following Taubes, every point $x$ in $\M^{\su}_{\pi}$ can be assigned an orientation 
$\epsilon(x)=\pm 1$ (take by convention $\epsilon([\mTheta])=+1$)
by considering the parity of the spectral flow of the \textit{Floer-Taubes
operator} $L^{\sw}$ from $[\mTheta]$ to $x$.  Fundamentally one would like to 
create an  topological invariant by taking 
the algebraic sum
\begin{equation}
    \sum_{x\in\M^{\su *}_{\pi}}\epsilon(x)
\end{equation}
as the Euler characteristic. However this sum can change with different choices 
of perturbation $\pi$ -- the phenomena of bifurcation of 
$\M^{\su}_{\pi}$ along 
$\M^{\su,I}_{\pi}$ (i.e. birth or death of points in $\M^{\su*}_{\pi}$). 
This is corrected by the addition of a 
counter-term associated to each point $x'$ in $\M^{\su,I}_{\pi}$. Along the
strata $\A^{I}_{F}$ the Floer-Taubes operator splits orthogonally as
\textit{normal} and \textit{tangential} components. (With some care 
this splitting can be made to extend over the trivial strata.) 
The \textit{normal operator} $N^{\su,I}$ is complex, the complex 
structure coming from the stabilizer of $\G_{F}$ along $\A^{I}_{F}$. 
$N^{\su,I}$ is essentially a twisted version of the 
signature operator in 3-dimensions. Denote 
by $\SF^{\su,I}_{\nu}([\mTheta],x)$ the complex spectral flow of $N^{\su,I}$ from
$[\mTheta]$ to $x$ in $\A^{0,I}_{F}/\G_{F}$. If this did not depend on 
the path chosen then the expression
\begin{equation}
    \sum_{x'\in\M^{\su,I}_{\pi}}\epsilon(x')\SF^{\su,I}_{\nu}([\mTheta],x')
\end{equation}
would be the sum total of the required counter-terms. 
Unfortunately the proposed counter-term does depend on the path chosen.
This dependency is traced back to the fact that
$\pi_{1}(\A^{0,I}_{F}/\G_{F})\cong\Z$ and the value of $\SF_{\nu}^{\su,I}$ around
the generator is $\pm 2$. This is the
\textit{(spectral flow) obstruction 
to defining an integer
$SU(3)$-Casson invariant utilizing only the non-degenerate (perturbed)
moduli space}.

Boden-Herald \cite{BH} solve this by firstly allowing only very small pertubations so as 
to compare the perturbed and unperturbed moduli space. This in turn 
enables using the Chern-Simons function on the unperturbed moduli space
to cancel the obstruction on the counter-terms. In this way a topological 
invariant is obtained, but it is no longer obviously an integer.
Boden-Herald-Kirk \cite{BHK} extract an integer invariant by a slight
modification of the preceding. However the definition still relies 
heavily on 
the unperturbed moduli space and is not very natural. 
Cappell-Lee-Miller \cite{CLM2} found an ingenious solution involving a certain 
differential in the Floer homology chain complex. Their solution has the advantage 
of allowing large perturbations, but the definition seems 
unnatural and invokes a much more complicated object, Floer theory.

\subsection{$\mathbf{U(2)}$-Seiberg-Witten and 
$\mathbf{SU(3)}$-Casson}

Given the mentioned difficulties in establishing a satisfactory 
integral valued $SU(3)$-Casson invariant we propose an entirely
different approach using the moduli space for the parallel
\textit{$U(2)$-Seiberg-Witten (SW) equations}. Our main thesis may be summarized as follows:

\textit{The obstruction in $SU(3)$-Casson can be made
to cancel against the obstruction in (twice) $U(2)$-Seiberg-Witten. 
The SU(3)-Casson and U(2)-Seiberg-Witten equations when taken
in conjunction yield an integral invariant of integral homology 
spheres involving only the non-degenerate (perturbed) moduli spaces,
path independent spectral flow 
counter-terms and Atiyah-Patodi-Singer spectral invariants
and is not limited by small perturbations.}

Let us briefly explain the idea which is detailed in the main body of 
the article. (A preliminary study of this was treated in \cite{lim3}.)
Let $E\to Y$ be a $U(2)$-bundle and $S\to Y$ the 
complex spinor bundle. Denote by $\CC_{E}$ the space of pairs 
consisting of connection and section of $S\otimes E$ (the tensor 
product is taken over $\C$). The 
$U(2)$-Seiberg-Witten solutions are the critical points of the 
Chern-Simons-Dirac function $\csd\colon\CC_{E}\to\R$. The quotient space 
$\CC_{E}/\G_{E}$ by the gauge 
transformations $\G_{E}$ has highest strata the
\textit{irreducibles} $\CC^{*}_{E}/\G_{E}$, and the lowest the trivial 
connections $\CC^{0}_{E}/\G_{E}$. There
are two intermediate strata relevent to us, we denote them as Type I 
and Type II. 

Type I reducibles are those where the spinor component 
is zero. Type II reducibles are essentially configurations where the 
connection gives a parallel reduction $E=L_{0}\oplus L_{1}$ and one of 
the components of the spinor in $S\otimes E=(S\otimes L_{0})\oplus 
(S\otimes L_{1})$ is identically zero.  

We have a corresponding decomposition of the Seiberg-Witten 
moduli space (as before we need to introduce non-degenerate perturbations $\pi'$):
\begin{equation}
    \M^{\sw}_{\pi'}=
    \M^{\sw *}_{\pi'}\cup\M^{\sw,I}_{\pi'}\cup\M^{\sw,II}_{\pi'}\cup\{[\mTheta]\}.
\end{equation}
As the metric and/or perturbation is varied bifurcation phenomena 
again happens along the lower stratas.
Since Type I reducibles are configurations where the spinor component is zero, we 
can identify $\M^{\sw,I}_{\pi'}$ as the moduli space of (perturbed) flat 
$U(2)$-connections. Notice  that 
$\A^{0,1}_{F}/\G_{F}=\CC^{0,1}_{E}/\G_{E}$.
It is straightforward to arrange the perturbations 
so that we have the identification
\begin{equation}
    \M^{\su,I}_{\pi}=\M^{\sw,I}_{\pi'}
\end{equation}
In the SW portion of the theory we again need to consider counter-terms associated to
a point $x$ in $\M^{\sw,I}_{\pi'}$. This involves the complex spectral flow 
$\SF_{\nu}^{\sw,I}([\mTheta],x)$ from $[\mTheta]$ to $x$   
of a corresponding normal operator $N^{\sw,I}$ which is a twisted Dirac 
operator. Again 
we have a \textit{spectral flow obstruction} -- around a generator
of $\pi_{1}(\CC^{0,I}_{E}/\G_{E})\cong\Z$ the value of $\SF_{\nu}^{\sw,I}$ 
is $\pm 1$. Now following a well established tradition in index theory 
we play-off $\SF_{\nu}^{\su,I}([\mTheta],x)$ (essentially a signature 
operator) against $\SF_{\nu}^{\sw,I}([\mTheta],x)$ (a Dirac operator) by
working with the sum (after taking signs into account)
\begin{equation}\label{intro-sf-equ}
   \epsilon(x)\Bigl\{(\SF_{\nu}^{\su,I}+2\SF_{\nu}^{\sw,I})([\mTheta],x)+4c(g,\pi')\Bigr\}.
    \end{equation}
The extra $4c(g,\pi')$ is an Atiyah-Patodi-Singer (APS) \cite{APS} spectral invariant term inserted to 
suppress spectral-flow at $[\mTheta]$ under variation of the metric. 
When $\pi'=0$ this takes the form
\begin{equation}
    c(g,0)=\xi+\frac{1}{8}\eta(B)
\end{equation}
where $\xi$ is the APS spectral invariant for the (untwisted) Dirac operator on $Y$ and
$\eta(B)$ that for the (untwisted) signature operator. By index theorems this sum 
is always an integer and reduces $\bmod2$ to the $\mu$-invariant for 
$Y$.
The sum (\ref{intro-sf-equ}) is now independent of the path chosen and should appear as the
correct counter-term in a proposed invariant. 
It is then clear that in the highest strata we should base a 
topological invariant on the following combination of SW and Casson theories:
\begin{equation}
    \sum_{x\in\M^{\su*}_{\pi}}\epsilon(x)+
    2\sum_{x'\in\M^{\sw*}_{\pi'}}\epsilon(x').
\end{equation}
This does not quite complete the invariant for as we vary the metric 
and/or perturbation, in the SW-case there is also bifurcation along (i) the trivial 
strata (ii) the Type II strata (this is essentially the $U(1)$-SW 
moduli space). Case (ii) can be straightforwardly  handled by spectral flow 
terms which do not suffer from anomalies.

Case (i) is different in that bifurcations along the trivial strata 
\textit{give birth or death into the Type II strata}. Thus what is 
needed is a \textit{counter-term for the counter-term associated to 
$\M^{\sw,II}_{\pi'}$.} We identify such an expression which turns out 
to be
\begin{equation}
    c(g,\pi')( c(g,\pi')-1).
\end{equation}
This completes all the counter-terms 
required to create a topological invariant.

We finish this introduction with a technical remark.
The traditional approach to 
perturbations is to perturb the Chern-Simons(-Dirac) function. In the 
SW-case there are problems to getting an adequate many for 
transversality -- the author has yet to find a satisfactory solution 
this way. Instead we perturb the gradient of the Chern-Simons(-Dirac) 
function directly (which is a much simpler procedure) in a way which preserves
most of the features we would expect from gradient perturbations. 

In \S\ref{basics} we discuss the basics of 
${U(2)}$-Seiberg-Witten and ${SU(3)}$-Casson theory. This includes
introducing the class of admissible perturbations, compactness of the 
moduli space and the splitting of 
slice spaces along the reducible stata. We also analyze 
the behaviour of an admissible perturbation near a reducible stata,
especially its \textit{normal linearization}. We conclude with showing 
how to construct adequately many admissible perturbations so as to 
obtain a non-degenerate moduli space.

\S\ref{definition} begins with a discussion of the Floer-Taubes 
operator which is the operator which allows us to define orientations 
as well as counter-terms.
We show that along the reducibles the normal component of this operator
is self-adjoint and thus the notion of spectral-flow is defined.
The definition of the invariant now proceeds after a  
discussion of the counter-terms.
\S\ref{proof-main}:
Proof of the main theorem.
\S\ref{sec-orient}
Treats orientation issues which are crucial to obtaining the correct 
signs for all the terms in the definition of the invariant.

\section{$\mathbf{U(2)}$-Seiberg-Witten and $\mathbf{SU(3)}$-Casson 
theory}\label{basics}

\textbf{Standing Convention} \textit{Throughout this article $Y$ will denote an oriented closed 
integral homology $3$-sphere (ZHS). $Y$ will also be assumed to
be have a fixed Riemannian metric $g$.} 

\subsection{The Equations}
Let $P\to Y$ be the unique spin-structure on $Y$ (up to
equivalence). In the (real) Clifford
bundle $CL(T^{*}Y)\cong CL(Y)$ the volume form $\omega_Y$ has the property
that $\omega^2_{Y}=1$. The action of $\omega_Y$ on $CL(Y)$ induces
a splitting into $\pm 1$ eigenbundles $CL^{+}\oplus CL^{-}$.
Both $CL^{+}$ and $CL^{-}$ are bundles of algebras over $Y$
with each fibre isomorphic, as an algebra, to the quaternions $\qu$ 
(see for instance \cite{law}).
Let $S\to Y$ be the complex spinor bundle on which $CL^{+}$ acts non-trivially.
This is a rank 2 complex Hermitian vector bundle. 

Fix $E\to Y$ a (trivial) $U(2)$-vector bundle, i.e. a rank 2 Hermitian complex 
vector bundle. Twist $S$ by forming the tensor
product (over $\C$) $S\otimes E$, a rank 4 complex vector bundle.
The Clifford action on $S$ naturally extends to $S\otimes E$ by
the rule $\alpha\cdot(\phi\otimes e)=(\alpha\cdot\phi)\otimes e$.
The action of $\mLambda^{2}\otimes\ad E$ on $S\otimes E$
defines a fibrewise bilinear form $\{\cdot\}_{0*}$ on $S\otimes E$ by the rule
\begin{equation}
\langle \alpha\cdot\phi,\psi\rangle = \langle \alpha,\{\psi\cdot\phi\}_{0*}\rangle.
\end{equation}

The \textit{$U(2)$-Seiberg-Witten Equation (in 3-dimensions)} is
the equation defined for a pair $(A,\mPhi)$ consisting of a connection
on $E$ and a spinor $\mPhi$ (i.e. section of $S\otimes E$). The equation reads:
\begin{equation}\label{sw*-equ}
    F_A -\{\mPhi\cdot\mPhi\}_{0*}=0,\quad D_{A}\mPhi=0,
\end{equation}
where $F_A$ is the curvature of $A$, and since $A$ is an 
$U(2)$-connection, $F_{A}$ is a section of ${\mLambda}^2\otimes \ad E$.
$D_{A}$ is the twisted Dirac operator on $S\otimes E$ and $\{\cdot \}_{0*}$
the quadratic form above.

Let $F\to Y$ be a fixed (trivial) $SU(3)$-vector bundle over $Y$, that is a rank 3 
Hermitian complex vector bundle with trivialized determinant.
In \textit{$SU(3)$-Casson theory} we are concerned with the flat connections 
on $F$, i.e. the solutions of the \textit{flat equation}
\begin{equation}
    F_{A}=0
\end{equation}
where $A$ is an $SU(3)$-connection on $F$.

\subsection{Configuration Spaces}
$\CC_{E}$ will denote the \textit{configuration space} of pairs
$(A,\mPhi)$ where $A$ is unitary connection on $E$ and $\mPhi$
a twisted spinor, i.e. a section of $S\otimes E$. 
In this article we shall be working in an $L^{2}_{2}$-Sobolev gauge 
theory. This means both $A$ and $\mPhi$
shall be of class $L^{2}_{2}$. The \textit{gauge automorphism group}
$\G_{E}$ will be the $L^{2}_{3}$-sections of $\Ad E$, the unitary bundle
automorphisms of $E$. The action of $\G_{E}$ on $\CC_{E}$ is from the 
right as $g\cdot(A,\mPhi) = (g(A), g^{-1}\mPhi)$ where the convention is that
$g(A)$ is the pull-back connection. 
The \textit{$U(2)$-SW moduli space} $\M^{\sw}$ is the solutions of 
(\ref{sw*-equ}) modulo gauge equivalence.

In a likewise manner $\A_{F}$ is the space of 
$L^{2}_{2}$-Sobolev $SU(3)$-connections
on $F$. The gauge group of $L^{2}_{3}$ automorphisms we denote by
$\G_{F}$. The \textit{$SU(3)$-flat moduli space} is denoted by 
$\M^{\su}$.

$\CC_{E}$ is an affine space modelled on $L^{2}_{2}(\mLambda^{1}\otimes\ad E)
\times L^{2}_{2}(S\otimes E)$. Here $\ad E$ denotes the bundle of Hermitian
skew endomorphisms of $E$. The tangent space to the identity of $\G_{E}$
is $L^{2}_{3}(\ad E)$ and the derivative at the identity of the gauge orbit
map $\G_{E}\to\CC_{E}$, $g\mapsto g\cdot(A,\mPhi)$ is given by
the operator 
\begin{equation}\label{gauge-deriv-sw}
\begin{array}{rcl}
\delta_{A,{\mPhi}}^0\colon L^2_3(\ad E) &\to& 
L^2_2({\mLambda}^1\otimes\ad E)\oplus
L^2_2(S\otimes E),\\
\delta^0_{A,\mPhi}(\gamma) &=& (d_A\gamma,-\gamma(\mPhi) ).
\end{array}
\end{equation}
The \textit{slice space} at $(A,\mPhi)$ is the $L^{2}$-orthogonal to the
image of $\delta^{0}_{A,\mPhi}$ and is denoted by $X_{A,\mPhi}$.

$\A_{F}$ is an affine space modelled on $L^{2}_{2}(\mLambda^{1}\otimes \ad F)$.
The tangent space at the identity to $\G_{F}$ is $L^{2}_{3}(\ad F)$
and the derivative at the identity of the gauge orbit
map $\G_{F}\to\A_{F}$, $g\mapsto g\cdot A$ is given by
the operator 
\begin{equation}\label{gauge-deriv-su}
d_{A}^0\colon L^2_3(\ad F)\to 
L^2_2({\mLambda}^1\otimes\ad F).
\end{equation}
The \textit{slice space} at $A$ is the $L^{2}$-orthogonal to the
image of $d^{0}_{A}$ and is denoted by $X_{A}$. 

\textbf{Remark} In $L^{2}_{2}$-gauge theory in 3-dimensions the 
connections and spinors are continuous objects. Since we are in the 
continous range for Sobolev theory, the SW and flat equations are 
well-defined as equations in $L^{2}_{1}$.  Details of Sobolev gauge 
theory can be found in for instance Freed-Uhlenbeck \cite{FU}.

\subsection{Reducibles}
We shall call $(A,\mPhi)\in\CC_{E}$ \textit{reducible} if the stabilizer of
$(A,\mPhi)$ is non-trivial, otherwise we call $(A,\mPhi)$
\textit{irreducible}. Geometrically a reduction happens in two
ways. In the first case $\mPhi=0$; then $\stab(A,\mPhi)$ is at least $U(1)$
(the gauge transformations which are  multiplication by a complex unit).
The second is when $A$ is reducible
as $A_{0}\oplus A_{1}$ in a parallel splitting 
$E=L_{0}\oplus L_{1}$ 
and $\mPhi$ is $A$-\textit{reducible} in the sense that 
$\mPhi=(\phi_{0},\phi_{1})\in L^{2}_{2}(S\otimes L_{0})\oplus 
L^{2}_{2}(S\otimes L_{1})$
with at least one of $\phi_{0,1}=0$. There are various reducible strata
with stabilizers $U(1)$, $U(1)\times U(1)$ and $U(2)$ however we shall only
be concerned with the following ones:
\begin{itemize}
	\item Type I:  $\mPhi=0$ and $A$ is irreducible as a connection on $E$. In this case
	$\stab(A,\mPhi)=\stab(A)$ under the action of $\G_{E}$ and this is easily seen to be just 
	those which are multiplication by a complex unit. Thus  $\stab(A)\cong U(1)$.
	\item Type II:  $A$ is reducible as $A_{0}\oplus A_{1}$ 
	and $\mPhi=(\phi_{0},0)\neq 0$. The stabilizer
	consists of the gauge transformations which have block diagonal form
	$$
	\left(
	\begin{array}{cc}
	1&0\\
	0& g
	\end{array}
	\right).
	$$
	Since the gauge automorphisms are sections of $\Ad E$ (which is fibrewise 
	$\cong U(2)$)
	we see that $g\in U(1)$. Thus in this case $\stab(A,\mPhi)\cong U(1)$ again.
	\item Trivial:  $A$ is a trivial connection $\mTheta$ and $\mPhi=0$. Here
	$\stab(A,0)=\stab(A)\cong U(2)$.
\end{itemize}
As general notation the irreducible portion of $\CC_{E}$ shall be denoted
by $\CC_{E}^{*}$ and the reducible
portion $\CC_{E}^{r}$. The Trivial, Type
I and II reducible stratas shall be denoted by $\CC^{0}_{E}$, 
$\CC^{I}_{E}$ and $\CC^{II}_{E}$ respectively. 
Note that within 
our definition $\CC^{0}_{E}$, 
$\CC^{I}_{E}$ and $\CC^{II}_{E}$ are mutually disjoint.

Occasionally it will be useful to specify the splitting $E=L_{0}\oplus 
L_{1}$ in a Type II reducible. We denote  $\CC^{II}(L_{0},L_{1})\subset 
\CC^{II}$ the subset with the given splitting and with the spinor 
component in $S\otimes L_{1}$ vanishing. Under the action of 
$\G_{E}$, $\CC^{II}(L_{0},L_{1})$ sweeps out $\CC^{II}$.

Finally let 
$\M^{\sw*}=\M^{\sw}\cap \CC^{*}_{E}/\G_{E}$ and 
$\M^{\sw,r}=\M^{\sw}\backslash\M^{\sw*}$ denote the irreducible and 
reducible portions of the moduli space.

\begin{lemma}
The only possible reducible SW-solutions on the ZHS $Y$ are of Type I, II or Trivial.
The Type I reducibles correspond to the irreducible solutions of the flat equation
$F_{A}=0$ on $E$. The Type II reducibles correspond to the solutions of the 
$U(1)$-SW-equation. 
\end{lemma}

\proof (sketch) If $\mPhi=0$ then the SW-equation reduces to the flat equation
$F_{A}=0$ and the only reducible solutions on a ZHS are trivial ones. 
The stabilizer of an irreducible solutions is clearly $U(1)$.
On the other hand if $A=A_{0}\oplus A_{1}$ is reducible and 
$\mPhi=(\phi_{0},\phi_{1})\in L^{2}_{2}(S\otimes L_{0}) \oplus L^{2}_{2}(S\otimes L_{1})$
is $A$-reducible with say $\phi_{1}=0$ then the $U(2)$-SW-equation reduces to the
two sets of equations:
(i) $F_{A_{0}}=\{\phi_{0}\cdot\phi_{0}\}_{0}$, $D_{A_{0}}\phi_{0}=0$. This is the
$U(1)$-SW-equation. Assume $\phi_{0}\neq 0$, otherwise as before on a
ZHS we are back in the trivial solution (ii) $F_{A_{1}}=0$, this clearly has only the
trivial solution. Thus, apart from Type I and trivial reducibles, we only get
Type II. \qed

A reducible $SU(3)$-connection admits a parallel splitting $A=A_{0}\oplus A_{1}$
corresponding to $F=F_{0}\oplus F_{1}$. We refer to a $U(2)$-\textit{reducible} 
or Type I as one for which  $F_{0}$ is an $U(2)$-bundle and $A_{0}$ is irreducible.
Since $F$ is an $SU(3)$-bundle, this forces $F_{1}$ to be 
$\cong\overline{{\rm det}F_{0}}$
and $A_{1}$ the connection induced by $A_{0}$. The stabilizer of a $U(2)$-reducible
can be verified to be $\cong U(1)$. If additionally $A_{1}$ is actually trivial we 
term $A$ to be $SU(2)$-\textit{reducible}. 

The strata of irreducibles 
is denoted $\A^{*}_{F}$, the
Type I reducibles
by $\A^{I}_{F}$ and the strata of trivial connections by $\A^{0}_{F}$.
Set $\M^{\su*}=\M^{\su}\cap\A^{*}_{F}/\G_{F}$ and 
$\M^{\su,r}=\M^{\su}\backslash\M^{\su*}$.

The following is clear:

\begin{lemma} 
On the ZHS $Y$ the only reducible solutions to the flat $SU(3)$-equations
are flat $U(2)$-reducibles (in fact $SU(2)$-reducible) and trivial connections.
\end{lemma}

\textbf{Convention} It will often be convenient to simultaneously 
treat both the SW and SU(3) theories. To this end we shall employ the 
notation $\ZZ$ for either $\CC_{E}$ or $\A_{F}$, and $\G$ the 
corresponding group of gauge transformations.

\subsection{Admissible Perturbations}

The gauge group $\G_{E}$ acts naturally on the tangent space
$L^2_2(\mLambda^1\otimes\ad E)\times L^2_2(S\otimes  E)$ by
conjugation in the fibres of $\ad E$ and directly on $E$.
Thus the notion of a $\G_{E}$-equivariant map $\CC_{E}\to
L^2_2(\mLambda^1\otimes\ad E)\times L^2_2(S\otimes  E)$ makes sense.
Define an {\it admissible perturbation} $\pi$ on $\CC_{E}$ to be a 
$C^{3}$ $\G_{E}$-equivariant map
$\pi=(*k,l)\colon\CC_{E}\to L^2_2(\mLambda^1\otimes\ad E)\times L^2_2(S\otimes  E)$ 
satisfying
\begin{enumerate}
\item
$\pi_{A,\mPhi}\in X_{A,\mPhi}$
\item
the linearization (i.e. derivative) $(L\pi)_{A,\mPhi}$ at $(A,\mPhi)$ is a bounded linear
operator from $L^{2}_{2}(\mLambda^1\otimes\ad E)\oplus L^2_2(S\otimes 
E)$ back to itself
\item\label{uniform-bd}
there is a uniform bound 
$$
\normi{\pi_{A,\mPhi}}{L^{2}_{2,A}}=
\sum^{2}_{i=0}\normi{(\nabla^{A})^{i} 
k_{A,\mPhi}}{L^2}+\normi{(\nabla^{A})^{i}l_{A,\mPhi}}{L^2}\le C
$$
\item $\pi$ has support contained in $\CC^{*}_{E}\cup \CC^{0,I,II}_{E}$.
\item $\pi$ depends only on the spinor component $\mPhi$ in a neighbourhood of 
the trivial orbit.
\end{enumerate}
If $\pi=(*k,l)$ as above then we perturb
the $U(2)$-Seiberg-Witten equation by setting
\begin{equation}\label{sw-eq2}
F_A -\{\mPhi\cdot\mPhi\}_{0*}+k_{A,\mPhi}=0,\quad
D_A {\mPhi} +l_{A,\mPhi}=0.
\end{equation}
The corresponding moduli space is denoted $\M^{\sw}_{\pi}$,
the irreducucible portion $\M^{\sw*}_{\pi}$, the reducible portion 
$\M^{\sw,r}_{\pi}$ etc.

\begin{lemma}
The only possible reducible perturbed $U(2)$-SW-solutions on the ZHS 
$Y$ are of Type I, II or Trivial.
The Type I reducibles correspond to the irreducible solutions of the 
perturbed flat equation
$F_{A}+k_{A}=0$ on $E$. The Type II reducibles correspond to the solutions of a 
perturbed $U(1)$-SW-equation. 
\end{lemma}

\proof The admissible perturbations have by definition support
in $\CC^{*}_{E}\cup \CC^{0,I,II}_{E}$ and therefore
no new kinds of reductions are introduced. \qed

An {\it admissible perturbation} $\pi'$ on $\A_{F}$ consists of a 
$C^{3}$ $\G$-equivariant map 
$\pi'\colon\A_{F}\to L^2_2(\mLambda^1\otimes\ad F)$ with $\pi'_{A}\in 
X_{A}$ and satisfying the parallel conditions stated above, i.e.
drop the spinor component $S\otimes E$. In particular $\pi'$ should 
be zero in a neighbourhood of the trivial orbit.
The perturbed flat $SU(3)$-equation now reads as
\begin{equation}
	F_{A} + *\pi_{A}' =0
	\label{pert-flat-eq}
\end{equation}
and the corresponding perturbed moduli spaces
$\M^{\su}_{\pi'}$, $\M^{\su*}_{\pi'}$, $\M^{\su,r}_{\pi'}$ etc.

\begin{lemma}
On the ZHS $Y$ the only reducible solutions to the perturbed flat $SU(3)$-equations
are perturbed flat $U(2)$-reducibles and trivial connections.
\end{lemma}

\begin{remark}\rm
The $SU(2)$-reducibles which happen in the unperturbed case
in general cease to remain so in the perturbed case. (i.e. become
$U(2)$-reducible but not $SU(2)$-reducible.)
\end{remark}

\subsection{Compactness}

Fix a smooth connection $\nabla^{0}$. A metric on $\CC_{E}/\G_{E}$ 
which induces the (quotient) topology is defined by the
rule
\begin{equation}
d([A,\mPhi],[A',\mPhi']) =\inf_{g\in\G_{E}}\Bigl\{ \sum^{2}_{i=0}
\normi{((\nabla^{0})^{i}(A-g(A')),\mPhi-g^{-1}\mPhi')}{L^{2}}\Bigr\}.
\end{equation}
Thus a subset $\mathcal{N}\subset\CC_{E}/\G_{E}$ is compact if and 
only if given any sequence $(A_{i},\mPhi_{i})$ such that the orbits
$[A_{i},\mPhi_{i}]\in\mathcal{N}$ there exists a subsequence 
$\{i'\}\subset\{i\}$ and gauge transformations $g_{i'}$ such that
$g_{i'}(A_{i'},\mPhi_{i'})$ converges in $L^{2}_{2}$.
In a similiar way a metric is defined on $\A_{F}/\G_{F}$.

\begin{prop}\label{compact-prop}
For any admissible perturbation $\M^{\sw}_{\pi}$ and 
$\M^{\su}_{\pi'}$ are compact subspaces.
\end{prop}
    
We shall not go through this in detail but refer to \cite{lim3} where
the proof applies in this context. (Property (\ref{uniform-bd}) in 
the definition of an admissible perturbation plays the crucial role.)

\subsection{The Fundamental Elliptic Complex}\label{sect-elliptic}

Following Taubes we should interpret the $U(2)$-SW-equation as the
zeros of the gauge equivariant \lq $L^{2}$-vector field\rq\ on $\CC_{E}$
\begin{equation}
{\cal X}^{\sw}(A,\mPhi) \stackrel{\rm def}{=} (*F_A -*\{\mPhi\cdot\mPhi\}_0, D_A\mPhi).
\end{equation}
This descends to the vector field $\widehat{\cal X}^{\sw}$ on 
$\CC^{*}_{E}/\G_{E}$ and the zeros are exactly $\M^{\sw*}$. 
Let ${\cal X}^{\sw}_{\pi}={\cal X}^{\sw}+\pi$, the perturbation of ${\cal X}^{\sw}$. 
The linearization of ${\cal X}^{\sw}_{\pi}$ at $x=(A,\mPhi)$ is
\begin{eqnarray}
    &\delta^{1,\pi}_{x}\colon L^{2}_{2}(\mLambda^1\otimes\ad E)\oplus L^2_2(S\otimes  E)
    \to
    L^{2}_{1}(\mLambda^1\otimes\ad E)\oplus L^2_1(S\otimes  E)\\
    &(a,\phi)\mapsto(*d_{A}a-\{\phi\cdot\Phi\}_{0}, 
    D_{A}\phi+a\cdot\mPhi)+(L\pi)_{A,\mPhi}(a,\phi).\nonumber
\end{eqnarray}
At a solution $x=(A,\mPhi)$ this fits into an (partial) elliptic complex:
\begin{equation}\label{ellip-seq1}
    L^{2}_{3}(\ad 
    E)\stackrel{\delta^{0}_{x}}{\longrightarrow}
    L^{2}_{2}(\mLambda^1\otimes\ad E)\oplus L^2_2(S\otimes  E)
    \stackrel{\delta^{1,\pi}_{x}}{\longrightarrow} 
    X_{x}\cap L^{2}_{1}.
\end{equation}
The (harmonic) cohomologies we denote by $\bH^{\sw,i}_{x}$, $i=0,1,2$. 
The \textit{non-degeneracy} (or \textit{regularity}) of $x$ is the condition 
$\bH^{\sw,2}_{x}=\{0\}$. By equivariance if $x$ is non-degenerate then
so are all points in the orbit of $x$ and so non-degeneracy of the 
orbit $[x]$ makes sense.

Let $\XX^{\su}$ denote the map $\A_{F}\to L^{2}_{1}(\mLambda^{1}\otimes\ad F)$,
$A\mapsto *F_{A}$. This is perturbed as $\XX^{\su}_{\pi'}=\XX^{\su}+\pi'$.
This has the linearization 
\begin{equation}
    *d^{1,\pi'}_{A}\colon L^{2}_{2}(\mLambda^{1}\otimes \ad F)
    \to  L^{2}_{1}(\mLambda^{1}\otimes \ad F),
    \quad
    a\mapsto *d_{A}a +(L\pi')_{A}(a).
\end{equation}
At a solution the elliptic complex in this instance is
\begin{equation}\label{ellip-seq2}
    L^{2}_{3}(\ad F)\stackrel{d^{0}_{A}}{\longrightarrow}
    L^{2}_{2}(\mLambda^{1}\otimes\ad F)\stackrel{*d^{1,\pi'}_{A}}{\longrightarrow}
    X_{A}\cap L^{2}_{1}
\end{equation}
with (harmonic) cohomologies $\bH^{\su,i}_{A}$, $i=0,1,2$.
Non-degeneracy is defined just as above.

\begin{lemma}
    $\dimn\bH^{\sw,1}_{x}=\dimn\bH^{\sw,2}_{x}$ and 
    $\dimn\bH^{\su,1}_{x}=\dimn\bH^{\su,2}_{x}$.
\end{lemma}

\proof
In the SW-case: the partial elliptic complex (\ref{ellip-seq1}) can be extended to a full 
one by replacing the last term with
\begin{equation}\label{full-1}
    \dots\stackrel{\delta^{1,\pi}_{x}}{\longrightarrow} 
    L^{2}_{1}(\mLambda^1\otimes\ad E)\oplus L^2_1(S\otimes  E)
    \stackrel{\delta^{0*}_{x}}{\longrightarrow} 
    L^{2}(\ad E)
\end{equation}
where $\delta^{0*}_{x}$ is the formal $L^{2}$-adjoint of $\delta^{0}_{x}$. 
This is an elliptic complex on an odd dimensional manifold and thus 
has zero index. Finally note that $\bH^{\sw,0}_{x}$ is the same as the
3rd-cohomology of the full complex. This proves the lemma in the SW-case. The same argument
holds in the $SU(3)$-Casson case; this time the full version of 
(\ref{ellip-seq2}) is extended by
\begin{equation}\label{full-2}
    \dots\stackrel{d^{1,\pi}_{x}}{\longrightarrow} 
    L^{2}_{1}(\mLambda^1\otimes\ad F)\stackrel{d^{0*}_{x}}{\longrightarrow} 
    L^{2}(\ad F).
\end{equation}
\qed

It follows from the 
Kuranishi local model that if $[x]$ is a non-degenerate point in $\M^{\sw}_{\pi}$ or
$\M^{\su}_{\pi'}$ then $[x]$ is an isolated point in the moduli space.

\subsection{Slice Splittings along Reducible Stratas}\label{splitting}

Let $x$ be a point in a reducible strata $\mathcal{R}$ in $\ZZ=\CC_{E}$ or $\A_{F}$.
Then there is an $L^{2}$-splitting of the slice space at $x$,
\begin{equation}
    X_{x}=X^{\tau}_{x}\oplus X^{\nu}_{x}
\end{equation}
into a \textit{tangential} component (superscripted $\tau$) and a 
\textit{normal}
component (superscripted $\nu$). The tangential component is 
essentially the slice space for the gauge action on $\mathcal{R}$
and this determines the normal component by taking the 
$L^{2}$-orthogonal. A precise way of describing this is as follows. At
$x$, $\stab(x)$ acts on $X_{x}$; the latter can be decomposed into 
an invariant subspace on which the action is trivial and an invariant
subspace on which the action is non-trivial. The first subspace is 
$X^{\tau}_{x}$ and the second $X^{\nu}_{x}$.

The splitting is actually induced at the level of the fibers of the
various vector bundles involved. In the $U(2)$-SW case the tangent 
space to the configuration space are the $L^{2}_{2}$-sections of
$V=(\mLambda^{1}\otimes\ad E) \oplus (S\otimes E)$.
$\G_{E}$ acts on each fiber of $V$ 
by conjugation $v\mapsto gvg^{-1}$ on the first factor and 
$\phi\mapsto g^{-1}\phi$ in the second. Then in the manner above, the 
invariant factors of the action of $\stab(x)$ give rise to a 
parallel splitting (with respect to the connection component of $x$) of the form 
\begin{equation}
V=V^{\tau}\oplus V^{\nu}
\end{equation}
where $V^{\tau}$ is the factor on which $\stab(x)$ acts trivially.
Then $X^{\tau}_{x}=X_{x}\cap L^{2}_{2}(V^{\tau})$ and $X^{\nu}_{x}=X_{x}\cap 
L^{2}_{2}(V^{\nu})$. In the $SU(3)$ case the relevant bundle 
$V=\mLambda^{1}\otimes\ad F$ and $\G_{F}$ acts on this by conjugation 
in each fiber.

Let $W$ be $\ad E$ in the SW-case and $\ad F$ in the $SU(3)$-case. The
action of $\stab(x)$ on $W$ in the same way as above also gives an 
$L^{2}$-decomposition
\begin{equation}
    W=W^{\tau}\oplus W^{\nu}.
\end{equation}
We may identify the components of the splittings in terms of the parallel splitting of 
$E$ or $F$ determined by the connection component of $x$.
For future reference we determine explicitly $W^{\tau}$, $W^{\nu}$, $V^{\tau}$ and $V^{\nu}$ 
in the cases that interest us. We precede this by a standard lemma.

\begin{lemma}
    (a) Suppose that $A$ is a reducible $U(2)$-connection on $E$ in a parallel splitting $E=L_{0}\oplus 
    L_{1}$. Then this induces a parallel splitting $\ad 
    E=i\R_{0}\oplus i\R_{1}\oplus(L_{0}\otimes\overline{L_{1}})$. 
    Here the $i\R_{j}$ factor 
    is the subbundle  $\ad E$ which is
    multiplication by pure imaginary constants on $L_{j}$.
    (b) Suppose that $A$ is a reducible $SU(3)$-connection on $F$ in a 
    parallel splitting $E\oplus L$ where $E$ is a $U(2)$-bundle. 
    Then $\ad F$ has a parallel splitting with respect to $A$ as 
    $\ad E\oplus (E\otimes \overline{L})$ where the
    $\ad E$ factor is the natural subbundle induced by the inclusion
    $E\subset F$.
    \end{lemma}

In the following table a vector space denotes 
the trivialized bundle over $Y$ with fibre that vector space.

\footnotesize

\begin{tabular}{|c|c|c|c|c|c|}\hline
&    Reducible& Parallel splitting &Adjoint bundle& $W^{\tau}$ & $W^{\nu}$ \\ 
\hline
  SW & Trivial & $E=\C^{2}$ & $\mathrm{u}(2)$
  &$\{0\}$ & $\mathrm{u}(2)$\\ \cline{2-6}
    &Type I & none &$\ad E$& $\ad E$ & $\{0\}$\\ 
    \cline{2-6}
    &Type II & $E=L_{0}\oplus L_{1}$ &
    $i\R_{0}\oplus i\R_{1}\oplus(L_{0}\otimes\overline{L}_{1})$
    & $i\R_{0}\oplus i\R_{1}$ &$L_{0}\otimes\overline{L}_{1}$\\ \hline
$SU(3)$ &Trivial & $F=\C^{3}$ &$\su (3)$
& $\{0\}$ & $\su (3)$\\ 
\cline{2-6}
&Type I &  $F=E\oplus L$ &$\ad E\oplus (E_{0}\otimes \overline{L})$&
$\ad E$ &$E\otimes\overline{L}$ \\ \hline
\end{tabular}

\begin{tabular}{|c|c|c|c|}\hline
&    Reducible& $V^{\tau}$ & $V^{\nu}$ \\ 
\hline
  SW & Trivial &  $\{0\}$ & 
  $(\mLambda^{1}\otimes\mathrm{u}(2))\oplus(S\otimes \C^{2})$\\ \cline{2-4}
    &Type I & $\mLambda^{1}\otimes\ad E$ & $S\otimes E$\\ 
    \cline{2-4}
    &Type II &  $(\mLambda^{1}\otimes i\R_{0})\oplus (S\otimes L_{0})$ &
$(\mLambda^{1}\otimes(L_{0}\otimes\overline{L}_{1}))\oplus(S\otimes L_{1})$\\ 
&& $\oplus (\mLambda^{1}\otimes i\R_{1})$&\\ \hline
$SU(3)$ &Trivial & $\{0\}$ & $\mLambda^{1}\otimes\su (3)$\\ 
\cline{2-4}
&Type I &   $\mLambda^{1}\otimes\ad E$ &
$\mLambda^{1}\otimes(E\otimes\overline{L})$ \\ \hline
\end{tabular}
\normalsize

\subsection{Normal Linearization of Perturbations}

Now that we have a splitting of the slice spaces along a reducible strata, 
we wish to establish a corresponding splitting for the linearization 
\begin{equation}\label{lin-split}
    (L\pi)_{x}=(L^{\nu}\pi)_{x}\oplus (L^{\tau}\pi)_{x}
\end{equation}
where the first factor maps into $L^{2}_{2}(V^{\nu})$ and the latter 
$L^{2}_{2}(V^{\tau})$. This and the crucial fact that $(L\pi)_{x}$
is symmetric are proven below, see Proposition~\ref{nor-sym-prop}.

\begin{lemma}\label{commute-lem}
Let $x\in\ZZ^{r}$, $g\in\stab (x)$ and
$\gamma\in \underline{\stab}(x)$, the Lie Algebra. Then for all $v$,
\begin{eqnarray*}
(L\pi)_{x}(g\cdot v) &=& g\cdot 
(L\pi)_{x}(v),\label{red-per-eq1}\\
(L\pi)_{x}(\gamma\cdot v) &=& \gamma\cdot 
(L\pi)_{x}(v).\label{red-per-eq2}
\end{eqnarray*}
\end{lemma}

\proof The condition that $\pi$ is equivariant with respect to the 
gauge action means $\pi_{g(x+sv)}=g\cdot\pi_{x+sv}$, $g\in\stab(x)$, 
$s\in\R$. Note that $g(x+sv)=x+sg\cdot v$ since $g\in\stab(x)$. Differentiating 
with respect to $s$ and evaluating at $s=0$ gives
\begin{equation}\label{red-per-eq4}
    (L\pi)_{x}(g\cdot v) = g\cdot(L\pi)_{x}(v).
    \end{equation}
The second relation of the lemma 
is obtained by varying $g$ in (\ref{red-per-eq4}).
\qed

\begin{lemma}\label{nor-sym-lem}
    Let $x\in\ZZ^{r}$. Then $(L\pi)_{x}$ is a 
    symmetric operator on $L^{2}_{2}(V)$, i.e.
    \begin{equation}
	\langle{(L\pi)_{x}(v),w}\rangle_{L^{2}} = 
	\langle{v,(L\pi)_{x}(w)}\rangle_{L^{2}},\quad v,w\in L^{2}_{2}(V).
    \end{equation}	
\end{lemma}

\proof
Let $\delta^{0}_{x}$ denote the linearized gauge action map
(\ref{gauge-deriv-sw}) or (\ref{gauge-deriv-su}) at $x$.
Let $v\in L^{2}_{2}(V)$. 
The definition of an admissible perturbation implies that
\begin{equation}\label{sym-lem-eq1}
    \langle{\pi_{x+sv},\delta^{0}_{x+sv}(\gamma)}\rangle_{L^{2}}=0
\end{equation}
for all $\gamma\in L^{2}_{3}(W)$. Observe that 
$\delta^{0}_{x+sv}(\gamma)=\delta^{0}_{x}(\gamma)+s\gamma\cdot v$.
Now let $\gamma\in\underline{\stab}(x)$ so that $\delta^{0}_{x}(\gamma)=0$.
Differentiating (\ref{sym-lem-eq1}) twice with respect to $s$ and 
putting $s=0$ gives the relation
\begin{equation}\label{sym-lem-eq2}
    \langle{(L\pi)_{x}(v),\gamma\cdot v)}\rangle_{L^{2}}=0.
    \end{equation}
If $x\in \CC^{I,II}_{E}$ or $\A_{F}^{I}$ then $\stab(x)\cong U(1)$ induces a 
complex structure on $V^{\nu}$; thus we can find
find a $\gamma$ such that $\gamma^{2}=-1$. Let 
$L=(L\pi)_{x}$. Applying (\ref{sym-lem-eq2}) to $\gamma(v)+w$
and invoking the skew-symmetry of $\gamma$ and Lemma~\ref{commute-lem}
gives
\begin{eqnarray*}
    0&=&\langle{L(\gamma(v)+w)),\gamma(\gamma(v)+w)}\rangle_{L^{2}}\\
    &=& \langle{L(v),w}\rangle_{L^{2}}-\langle{L(w),v}\rangle_{L^{2}}.
\end{eqnarray*}
If $x=\mTheta\in\A^{r}_{F}$ the definition of 
$\pi$ admissible requires $\pi$ to vanish near $\mTheta$ thus 
$(L\pi)_{\mTheta}=0$. If $x=(\mTheta,0)\in\CC^{r}_{E}$ the 
condition of being admissible requires $\pi_{A,\mPhi}$ near 
$(\mTheta,0)$ to depend only on $\mPhi$. Thus $(L\pi)_{\mTheta}$ 
acts only on the spinor factor $L^{2}_{2}(S\otimes E)$ of $V^{\nu}$ and this as above gets a 
complex struture from the action of the stabilizer. Repeating the 
argument we see that $(L\pi)_{\mTheta}$ is also symmetric. \qed

\begin{prop}\label{nor-sym-prop}
    Let $x\in \ZZ^{r}$. Then there exists an 
    $L^{2}$-orthogonal splitting 
    \begin{displaymath}
	(L\pi)_{x}=(L^{\tau}\pi)_{x}\oplus (L^{\nu}\pi)_{x}\colon
	L^{2}_{2}(V^{\tau})\oplus L^{2}_{2}(V^{\nu})\to
	L^{2}_{2}(V^{\tau})\oplus L^{2}_{2}(V^{\nu}).
	\end{displaymath}
    Furthermore the 
    {\rm normal linearization} $(L^{\nu}\pi)_{x}$ has the following 
    properties (i) it commutes with the action of $\stab(x)$ (ii) it is 
    symmetric with respect to the $L^{2}$-inner product
    (iii) if additionally $\pi=0$ on $\ZZ^{r}$ then  
    $(L^{\tau}\pi)_{x}=0$.
\end{prop}

Thus the \textit{normal linearization} of a perturbation is symmetric and complex linear in 
the case $\stab(x)\cong U(1)$ and symmetric and quaternionic linear 
when  $\stab(x)\cong U(2)$ corresponding to $x=\mTheta$ in the SW-case (and zero for $x=\mTheta$ in the 
$SU(3)$-case).

\proof
Firstly if $x\in \ZZ^{r}$ and $g\in\stab(x)$ then
$\pi_{x}=\pi_{g(x)}=g\cdot\pi_{x}$ implies that $\pi_{x}\in 
L^{2}_{2}(V^{\tau})$. Therefore for $w\in L^{2}_{2}(V^{\tau})$ we have
$(L\pi)_{x}(w)\in L^{2}_{2}(V^{\tau})$ as well. Now suppose $v\in 
L^{2}_{2}(V^{\nu})$. Then Lemma~\ref{nor-sym-lem} shows that
$\langle{(L\pi)_{x}(v),w}\rangle_{L^{2}}=\langle{v,(L\pi)_{x}(w)}\rangle_{L^{2}}=0$,
$w\in L^{2}_{2}(V^{\tau})$.
Since this holds for all $w\in L^{2}_{2}(V^{\tau})$ we deduce that 
$(L\pi)_{x}(v)\in L^{2}_{2}(V^{\nu})$. This proves the splitting. 
Items (i) and (ii) follow immediately from
Lemmas~\ref{commute-lem} and \ref{nor-sym-lem}.
Item (iii) is clear. \qed

\begin{prop}\label{ellip-split}
    The  complexes (\ref{ellip-seq1}) and (\ref{ellip-seq2}) at a 
    solution $x$ decompose orthogonally into {\rm tangential} 
    and {\rm normal} complexes:
\begin{eqnarray*}
L^{2}_{3}(W^{\tau})
\stackrel{\delta^{0,\tau}_{x}}{\longrightarrow} L^{2}_{2}(V^{\tau})
    \stackrel{\delta^{1,\tau,\pi}_{x}}{\longrightarrow}
    X^{\tau}_{x}\cap L^{2}_{1},\\
L^{2}_{3}(W^{\nu})
\stackrel{\delta^{0,\nu}_{x}}{\longrightarrow} L^{2}_{2}(V^{\nu})
    \stackrel{\delta^{1,\nu,\pi}_{x}}{\longrightarrow}
    X^{\nu}_{x}\cap L^{2}_{1}.
  \end{eqnarray*}
We have corresponding $L^{2}$-orthogonal splittings of cohomologies
as $\bH^{i}_{x}=\bH^{i,\tau}_{x}\oplus\bH^{i,\nu}_{x}$.
\end{prop}

\proof Follows directly from the preceding. \qed

\subsection{Abundance of Perturbations}

\begin{prop}\label{adm-lem-sw}
    There exists non-degenerate admissible perturbations, i.e.
$\pi$ ($\pi')$ such that
$\M^{\sw}_{\pi}$ ($\M^{\su}_{\pi'}$) consists entirely of 
non-degenerate points. Furthermore $\pi$ ($\pi'$) may be chosen
to have support in any arbitarily small gauge invariant neighbourhood of the
subspace of unperturbed SW ($SU(3)$-flat) solutions.
\end{prop}

As before ${\cal Z}$ denotes either $\CC_{E}$ or $\A_{F}$, and $\G$ the 
gauge group. The strategy is to construct perturbations locally in $\ZZ/\G$. Done correctly 
these will be admissible. To do so we  require some preliminary 
technical lemmas.
Introduce the notation
$B(\epsilon)$ for the $\epsilon$-ball in the slice space
$X_{x}$. Denote by $\beta\colon X_{x}\to[0,1]$ a smooth cut-off function
with support in $B(\epsilon)$. Let $\delta^{0}_{x}$ the zeroth differential in the elliptic 
complex and $W\to Y$ the adjoint bundle.

\begin{lemma}\label{implicit}
Fix $x\in{\cal Z}$. For all $\epsilon>0$ sufficiently small there is a differentiable function
$\xi\colon B(\epsilon)\times X_{x}\to (\kernel \delta^{0}_{x})^{\perp}
\subset L^{2}_{3}(W)$ such that given any 
$(\alpha,v)\in B(\epsilon)\times X_{x}$, the equation
\begin{equation}\label{xi-equ}
v +\delta^{0}_{x}\circ\xi (\alpha,v) \in X_{x+\alpha}
\end{equation}
holds. Here $(\kernel \delta^{0}_{x})^{\perp}$ denotes the 
$L^{2}$-orthogonal complement.
\end{lemma}

\proof Apply the Implicit Function theorem 
to the map
\begin{displaymath}
H(\xi,\alpha,v) =\delta^{0*}_{x+\alpha}(\delta^{0}_{x}(\xi)
+v)
\end{displaymath}
from $(\kernel \delta^{0}_{x})^{\perp}\times B(\epsilon)\times X_{x}
\to (\kernel \delta^{0}_{x})^{\perp}\cap L^{2}_{1}$. 
The linearization of $H$ at the origin  restricted to 
$(\kernel \delta^{0}_{x})^{\perp}$ is an isomorphism. This establishes
the existence of the function $\xi=\xi(\alpha,v)$ as claimed but only for $\alpha$ and $v$
defined in sufficiently small neighbourhoods of zero. However notice that
if $v$ satisfies (\ref{xi-equ}) then 
for any real constant $c$,
$cv$ satisfies the same equation
but with $\xi$ replaced by $c\xi$. That is we can allow the $v$
in $\xi$ to be defined for all $X_{x}$ by extending $\xi$ linearly in that factor.
\qed

\begin{lemma}\label{pert1-lem}
    Assume the hypotheses of Lemma~\ref{implicit}. Then for all
$\epsilon>0$ sufficiently small there is a constant $c$ (independent of 
$\alpha$, $v$ and $\epsilon$) such that $\normi{\xi(\alpha,v)}{L^{2}_{3}}\le 
c\normi{v}{L^{2}_{2}}$.
\end{lemma}

\proof
$\xi$ satisfies $H(\xi,\alpha,v)=0$.
Thus
\begin{equation}\label{D-equ}
\Delta_{x}\xi +N_{1}(\alpha,\xi) + N_{2}(\alpha,v) + 
\delta^{0*}_{x}(v)=0
\end{equation}
where $\Delta_{x}$ is the Laplacian $\delta^{0*}_{x}\delta^{0}_{x}$
and $N_{1}$ and $N_{2}$ are lower order terms. $N_{1}$ is a bilinear expression
in $\alpha$ and $\delta^{0}_{x}(\xi)$. $N_{2}$ is a bilinear expression
in $\alpha$ and $v$.
After some calculation it is seen that
$N_{1}$, $N_{2}$ satisfy, by Sobolev theorems
\begin{eqnarray}
\normi{N_{1}(\alpha,\xi)}{L^{2}_{1}}&\leq& \const \normi{\alpha}{L^{2}_{2}}
\normi{\xi}{L^{2}_{3}}\label{N-equ}\\
\normi{N_{2}(\alpha,v)}{L^{2}_{1}} &\leq &\const
\normi{\alpha}{L^{2}_{2}}\normi{v}{L^{2}_{2}}.\nonumber
\end{eqnarray}
On the other hand since $\Delta_{x}$ is invertible 
on $(\kernel\delta^{0}_{x})^{\perp}$,
\begin{equation}\label{D2-equ}
\normi{\xi}{L^{2}_{3}} \leq\const\normi{\Delta_{x}\xi}{L^{2}_{1}}.
\end{equation}
Now make $\epsilon>0$ sufficiently
small so that $\normi{\alpha}{L^{2}_{2}}$ is correspondingly small.
Then
(\ref{D-equ}), (\ref{N-equ}) and (\ref{D2-equ}) give
$\normi{\xi}{L^{2}_{3}} \leq\const\normi{v}{L^{2}_{2}}$. \qed

\begin{prop}\label{pert-1}
Assume $x\in{\cal Z}^{*}$.
Given any $v\in X_{x}$ there is an admissible perturbation $\pi$
such that $\pi_{x}=v$. Furthermore the support of $\pi$
may be chosen to be contained in an arbitarily small $\G$-invariant 
neighbourhood of the orbit  $\G\cdot x$.
\end{prop}

\proof
Identify the slice space $X_{x}$ with the actual slice $x+X_{x}$.
Let $\epsilon>0$ be sufficiently small so that 
$B(\epsilon)$ injects into the quotient space $\ZZ^{*}/\G$ and
the conclusions of 
Lemmas~\ref{implicit}, \ref{pert1-lem} hold. Construct a 
perturbation $\pi$ in $X_{x}$ by the rule that
\begin{equation}\label{pi-equ}
    \pi_{x+\alpha} = \beta(\alpha)v
+\delta^{0}_{x}\circ\xi(\alpha,\beta(\alpha)v),\quad \alpha\in 
B(\epsilon)\subset X_{x}.
\end{equation}
By construction $\pi$ has support in $B(\epsilon)$ and $\pi_{x+\alpha}\in 
X_{x+\alpha}$.
Extend $\pi$ to
$\ZZ^{*}$ by $\G$-equivariance. What remains is to show that $\pi$ is 
admissible provided $\epsilon$ is sufficiently small.
Equation (\ref{pi-equ}) and Lemmas~\ref{implicit}, \ref{pert1-lem} 
give a uniform bound
\begin{displaymath}
\normi{\pi_{x+\alpha}}{L^{2}_{2}}
\leq \const\normi{v}{L^{2}_{2}}\leq C.
\end{displaymath}
Here the Sobolev norm is taken with respect to some fixed connection
$A_{0}$, which is commensurate to the Sobolev norm taken to the 
connection component $A$ of $x$. Let $a$ denote the non-spinor 
component of $\alpha$. 
If $\normi{\alpha}{L^{2}_{2}}$ is sufficiently small then 
\begin{eqnarray*}
\normi{\nabla^{A+a}\pi_{x+\alpha}}{L^{2}}&\leq& \const
\normi{\nabla^{A}\pi_{x+\alpha}}{L^{2}},\\
\normi{\nabla^{A+a}\nabla^{A+\alpha}\pi_{x+a}}{L^{2}}&\leq& \const
\normi{\nabla^{A}\nabla^{A}\pi_{x+\alpha}}{L^{2}}
\end{eqnarray*}
uniformly. By reducing
$\epsilon$ again if necessary,
the uniform bound 
$\normi{\pi_{x'}}{L^{2}_{2,A'}} \le C$ where $A'$ is the connection 
component of $x'$, is established. \qed

\begin{prop}\label{pert-1a}
Assume $x\in{\cal Z}^{\alpha}$, $\alpha\in\{I,II\}$.
Given any $v\in X^{\tau}_{x}$ there is an admissible perturbation $\pi$
such that $\pi_{x}=v$. Furthermore the support of $\pi$
may be chosen to be contained in an arbitarily small $\G$-invariant 
neighbourhood of the orbit  $\G\cdot x$.
\end{prop}

\proof same argument as Proposition~\ref{pert-1}. \qed

Let $V$ be the vector underlying the configuration space $\ZZ$ 
(\S\ref{splitting}). We saw in Proposition~\ref{nor-sym-prop} that 
the normal linearization of any admissible perturbation 
defines a $\stab(x)$-equivariant 
symmetric bounded linear map $T\colon X^{\nu}_{x}\to 
X^{\nu}_{x}\subset L^{2}_{2}(V^{\nu})$. A certain 
assumption was placed on the definition of an admissible perturbation in order for 
Proposition~\ref{nor-sym-prop} to be valid: namely that $\pi$ did not 
depend on the connection component in a neighbourhood of the trivial 
connection orbit. This forces, when $x=\mTheta$ that $T$ is trivial 
on the first factor of $X^{\nu}_{\mTheta}=L^{2}_{2}(\mLambda^{1}\otimes\ad 
E)\oplus L^{2}_{2}(S\otimes E)$ in the SW-case, and completely trivial 
on $X^{\nu}_{\mTheta}$ in the SU(3)-case.
We shall call such $\stab(x)$-equivariant 
symmetric bounded linear maps $T\colon X^{\nu}_{x}\to X^{\nu}_{x}$ 
\textit{admissible}.

\begin{prop}\label{pert-2}
    Let $x$ be in a reducible strata $\ZZ^{r}$.
    Given any admissible map $T\colon X^{\nu}_{x}\to X^{\nu}_{x}$ there 
    exists an admissible perturbation $\pi$ such that $\pi=0$ on the
    reducibles and 
    $(L\pi)_{x}=(L^{\nu}\pi)_{x}=T$. Furthermore $\pi$ may be assumed to be 
    supported in an arbitarily small 
    $\G$-invariant neighbourhood of $x$.
\end{prop}

\proof
In the slice space $X_{x}$ define 
\begin{equation}
\pi_{x+\alpha} = \beta(\alpha)T(\alpha)
+\delta^{0}_{x}\circ\xi(\alpha,\beta(\alpha)T(\alpha)),\quad \alpha\in 
B(\epsilon)\subset X_{x}.
\end{equation}
As in the proof of Proposition~\ref{pert-1}
this defines an admissible perturbation with the desired 
properties provided $\epsilon$ is sufficiently small. \qed

\proofof{Proposition~\ref{adm-lem-sw}}
Consider first the SW-case. We proceed inductively up the various 
strata beginning with the trivial strata $\{[\mTheta]\}$. Here the normal 
operator $N_{\mTheta}=D\oplus D$ acts on $X^{\nu}_{\mTheta}=S\oplus S$.
Let $\Sigma$ be the unit $L^{2}$-sphere within $L^{2}_{2}(S)$. Denote 
by $\Pi_{0}$ the vector space of admissible perturbations which vanish on $\CC^{r}_{E}$. 
If we introduce a perturbation $\pi\in\Pi_{0}$ then  the corresponding 
normal operator $N^{\pi}_{\mTheta}=D^{\pi}\oplus D^{\pi}$ where 
$D^{\pi}=D+(L^{\nu}\pi)_{\mTheta}$. Denote by $\V\to \Sigma\times\Pi_{0}$
the vector bundle whose fibre at $(v,\pi)$ is the real 
$L^{2}$-orthogonal $\langle{iv}\rangle^{\perp}\subset 
L^{2}_{1}(S)$. Then $f(v,\pi)=D^{\pi}(v)$ is a section of this vector 
bundle. We claim that $f$ is a submersion along $f^{-1}(0,0)\cap 
(\Sigma\times\{0\})$. To see this, consider the derivative of $f$ at 
$(v,0)\in f^{-1}(0,0)$ in the direction $\pi\in\Pi_{0}$. 
We have $(Lf)_{v,0}(\pi)=(L^{\nu}\pi)_{\mTheta}$. By varying $\pi$ 
and invoking 
Proposition~\ref{pert-2} we see that $(Lf)_{v,0}$  must be 
surjective. Futhermore we can make $f$ a submersion by restricting 
$\pi$ to 
some finite dimensional subspace ${\cal H}\subset \Pi_{0}$.
By the Sard-Smale theorem there must exist a $\pi'\in {\cal H}$ (which we can 
assume arbitarily small) such that $f^{-1}(0,\pi')$ is cut out 
transversely in $\Sigma$. However if $(v,\pi')\in f^{-1}(0,\pi')$ then the 
symmetry and complex linearity of $D^{\pi'}$ forces both $v$ and $iv$ 
to be orthogonal to its image. This means $L(f|_{\Sigma})$ is not a 
submersion along $f^{-1}(0,\pi')$ which is a contradiction. Therefore
$f^{-1}(0,\pi')$ is empty and $D^{\pi'}$ is invertible. So for $\pi'$,
$\bH^{1}_{\mTheta}=\bH^{1,\nu}_{\mTheta}=\ker N^{\pi}_{\mTheta}=\{0\}$ 
and $\{[\mTheta]\}$ is a non-degenerate  and hence an isolated 
point in $\M^{\sw}_{\pi'}$. 

Next, invoking Proposition~\ref{pert-1a} we can find a perturbation 
(also denoted as $\pi'$) so 
that $\M^{I,II}_{\pi'}$ is non-degenerate \textit{within} $\CC^{I,II}_{E}$, i.e. the 
tangential cohomologies $\bH^{1,\tau}_{x}$ are trivial. To get the 
normal cohomologies $\bH^{1,\nu}_{x}$ to vanish as well repeat the 
argument as for the trivial strata -- the usage of perturbations 
vanishing on $\CC^{r}_{E}$ ensure that $\M^{I,II}_{\pi'}$ remains 
unchanged. At this stage $\M^{\sw,0,I,II}_{\pi'}$ is non-degenerate and isolated 
within $\M^{\sw}_{\pi'}$. Non-degeneracy for $\M^{\sw*}_{\pi'}$
is achieved by another perturbation of the sort in Proposition~\ref{pert-1}.

In the $SU(3)$-flat case, we proceed just as in the preceding except 
that we may skip the initial step of perturbing near the trivial 
connection. This is because $\bH^{1}_{\mTheta}\cong H^{1}(Y)$ is always trivial.
\qed

\section{Definition of the Invariant}\label{definition}

\subsection{The Floer-Taubes operator}\label{Floer-Taubes}

Let $x$ be a point in the configuration spaces $\CC_{E}$ or $\A_{F}$. As in the preceding 
subsection let $\delta^{0}_{x}$ denote the zeroth differential in the 
elliptic complex
and $V$ the vector bundle whose $L^{2}_{2}$-sections models the 
tangent space to the configuration space. Let $W$ denote $\ad E$ in 
the SW-case and $\ad F$ in the $SU(3)$-case.

The Floer-Taubes operator (at $x$) is the \lq roll-up\rq\ of
of the full version of the fundamental elliptic complex of \S\ref{sect-elliptic}.
It is the bounded operator from the
$L^{2}_{2}$-sections to the $L^{2}_{1}$-sections of $W\oplus V$ given
in block diagonal form:
\begin{equation}
L^{\pi}_{x} =
\left(
\begin{array}{cc}
0& \delta^{0 *}_{x} \\
\delta^{0}_{x}& \delta^{1,\pi}_{x}
\end{array}
\right).
\end{equation}
As before $\delta^{0 *}_{x}$ is the formal $L^{2}$-adjoint of 
$\delta^{0}_{x}$. By construction the kernel of $L^{\pi}_{x}$ is 
just $\bH^{0}_{x}\oplus \bH^{1}_{x}$ and the cokernel which we identify 
with the $L^{2}$-orthogonal of the image, $\bH^{0}_{x}\oplus \bH^{2}_{x}$.

The proofs in the next two subsections are below.

\subsubsection{Orientability}

In order to orient the moduli space we consider the determinant line
$\detind L^{\pi}$ of the family $L^{\pi}_{x}$ parameterized
by $x$. This determinant line is equivariant with
respect to the gauge action  and descends to a line bundle
denoted as $\detind\widehat{L}^{\pi}$ on the quotient space.
The \lq orientability\rq\ of the quotient space is a consequence of the 
following:

\begin{lemma}\label{det-orient-lem}
The line bundle $\detind \widehat{L}^{\pi}$ is orientable, i.e. the pull-back 
over any closed loop is a trivial line bundle over $S^{1}$.
\end{lemma}

With this lemma
it is possible to define, as in the manner of Taubes, a relative
sign between non-degenerate zeros of $\widehat{\XX}_{\pi}$
as the basis for a Poincare-Hopf index.

Let us now consider the Floer-Taubes operator along a reducible 
stratum $\mathcal{R}$. Let $x\in\mathcal{R}$. Then according to
\S\ref{splitting} we have a splitting of the bundles $W$ 
and $V$ into tangential and normal components. 
This induces a 
splitting of the Floer-Taubes operator
\begin{equation}
    L^{\pi}_{x} = K^{\pi}_{x}\oplus N^{\pi}_{x}
    \end{equation}
into \textit{tangential} and \textit{normal} components acting on 
sections of $W^{\tau}\oplus V^{\tau}$ and $W^{\nu}\oplus V^{\nu}$ 
respectively (details below). Note that if we denote by $K_{x}$ and $N_{x}$ the
operators in the unperturbed situation then 
$K^{\pi}_{x}=K_{x}+(L^{\tau}\pi)_{x}$ and 
$N^{\pi}_{x}=N_{x}+(L^{\nu}\pi)_{x}$ where we split $(L\pi)_{x}$ 
according to Proposition~\ref{nor-sym-prop}.
Furthermore $N^{\pi}_{x}$ commutes 
with the action of $\stab(x)$. 

\begin{lemma}
    The family of operators $\{N^{\nu}_{x}\}$ always extends continuously over the trivial strata
except in the SW-Type II case. In the SW-Type II case this becomes 
true after restrict the family to the subset of reductions 
$(A,\mPhi)=(A_{0}\oplus A_{1},(\phi,0))$ in a fixed splitting 
$E=L_{0}\oplus L_{1}$.
\end{lemma}

This follows from the explicit descriptions of the normal operators, 
below.

Henceforth we shall assume a fixed splitting $E=L_{0}\oplus L_{1}$ for 
all SW-Type II reductions, i.e. restrict to $\CC^{II}(L_{0},L_{1})$. 
This is without any loss of generality as any 
other such splitting can be moved to the reference one by a gauge 
transformation.

The orientability of the quotient space of $\mathcal{R}$
is claimed by the next lemma.

\begin{lemma}\label{det-red-lem}
    The determinant line $\detind K^{\pi}$ descends to an orientable line 
    bundle $\detind\widehat{K}^{\pi}$ over the quotient space of 
    $\mathcal{R}$.
\end{lemma}

\subsubsection{Spectral Flow}

The normal operator is a formally self-adjoint Fredholm operator 
when $\pi=0$. In particular it is a ($L^{2}$-)symmetric operator with domain 
the $L^{2}_{2}$-sections of $W^{\nu}\oplus V^{\nu}$. 
In the presence of a non-trivial perturbation term 
Lemma~\ref{nor-sym-prop}
asserts that the normal operator continues to be symmetric on 
$L^{2}_{2}$-sections.

The next proposition is an observation in \cite{taubes}.

\begin{prop}\label{spec-prop}
    The normal operator $N^{\pi}_{x}$ regarded as an unbounded operator on 
$L^{2}(W^{\nu}\oplus V^{\nu})$ is essentially self-adjoint. 
It has only a real discrete spectrum with no accumulation points. Each eigenvalue is 
of finite multiplicity and the eigenvalues are unbounded in both 
directions in $\R$. The normal operator depend differentiably on the
parameters $x$ and $\pi$.
\end{prop}

An immediate consequence of this proposition is that the concept of 
spectral flow is well-defined for $N^{\nu}_{x}$. The spectral flow for 
a path $\gamma$ in the
quotient space is taken to be the spectral-flow along any lift of 
$\gamma$; the independence of lift is clear because the operators
$N^{\pi}_{x}$ and $N^{\pi}_{g(x)}$ are conjugate to each other.

Let $\SF^{\sw,I}_{\nu}$, $\SF^{\sw,II}_{\nu}$, $\SF^{\su,I}_{\nu}$ denote the complex
spectral-flow in the various cases indicated.

\textbf{Convention} When working with spectral flow the initial and final operators 
in the family may have non-trivial kernel. Our convention will be the 
spectral flow across $-\epsilon^{2}$ where $0 <\epsilon\ll 1$.

\begin{lemma}\label{sf-lem}
{\rm (a)} Let  $\gamma(t)$, $t\in[0,1]$ be a  piecewise
differentiable loop in $\CC^{0,I}_{E}/\G_{E}=\A^{0,I}_{F}/\G_{F}$. 
Then
\begin{displaymath}
\SF^{\su,I}_{\nu}(\gamma) = 2\cs(\hat{\gamma}(0))-2\cs(\hat{\gamma}(1))=
-2\SF^{\sw,I}_{\nu}(\gamma).
\end{displaymath}
Here $\hat{\gamma}$ is a lift of $\gamma$ and
$\cs$ is the Chern-Simons function on $\CC^{0,I}_{E}=
\A^{0,I}_{F}/\G_{F}$
defined with respect to some fixed trivial basepoint connection.
\newline
{\rm (b)} If $\gamma$ is a loop in $\CC^{0,II}_{E}/\G_{E}$
then $\SF^{\sw,II}_{\nu}(\gamma)=0$.
\end{lemma}

\begin{remark}\rm
It is known that
$\pi_{1}(\A^{0,I}_{F}/\G_{F})\cong\Z$ and if $\gamma$ is 
a generator then $\cs(\hat{\gamma}(0))-\cs(\hat{\gamma}(1))=\pm 1$.
This gives us our spectral flow calculation $\SF^{\su,I}_{\nu}(\gamma)=\pm 
2$ and $\SF^{\sw,I}_{\nu}(\gamma)=\pm 1$ for a generator $\gamma$.
\end{remark}

\subsubsection{Details of the splittings and proofs}\label{detailproof}

Let us first summarize the splittings stated in \S~\ref{splitting}, in a more 
convenient form.

\footnotesize

\begin{tabular}{|c|c|c|c|}\hline
&    Reducible& $W^{\tau}\oplus V^{\tau}$ & $W^{\nu}\oplus V^{\nu}$ \\ 
\hline
  SW & Trivial &  $\{0\}$ & $(\mLambda^{0+1}\otimes{\rm u} (2))\oplus(S\otimes 
  \C^{2})$\\ \cline{2-4}
    &Type I & $\mLambda^{0+1}\otimes\ad E$ & $S\otimes E$\\ 
    \cline{2-4}
    &Type II &  $(\mLambda^{0+1}\otimes i\R_{0})\oplus 
(S\otimes L_{0})$ &
$(\mLambda^{0+1}\otimes(L_{0}\otimes\overline{L}_{1}))\oplus(S\otimes L_{1})$\\ 
&& $\oplus (\mLambda^{0+1}\otimes i\R_{1})$&\\ \hline
$SU(3)$ &Trivial & $\{0\}$ & $\mLambda^{0+1}\otimes\su (3)$\\ 
\cline{2-4}
&Type I &   $\mLambda^{0+1}\otimes\ad E$ &
$\mLambda^{0+1}\otimes(E\otimes\overline{L})$ \\ \hline
\end{tabular}

\normalsize

The corresponding tangential operators in the unperturbed situation is 
as below:

\footnotesize
\begin{tabular}{|c|c|c|}\hline
&    Reducible& Tangential Operator \\ 
\hline
  SW &Type I & $K^{\sw,I}_{A}(\gamma,a)=(d^{*}_{A},*d_{A}a+d_{A}\gamma)$\\ \cline{2-3}
    &Type II &  $K^{\sw,II}_{\alpha,\psi}(\gamma,a,\phi,\xi,b)=
    (d^{*}_{\alpha}a-B(\phi,\psi),*d_{\alpha}a-\{\phi\cdot\psi\}_{0*}+
    d_{\alpha}\gamma,$\\ 
&&$ D_{\alpha}\phi+a\cdot\psi-\gamma\psi,d^{*}b,*db+d\xi)$ \\ \hline
$SU(3)$ &Type I &   $K^{\su,I}_{\alpha}(\gamma,a)=(d^{*}_{\alpha}a,
*d_{\alpha}a+d_{\alpha}\gamma)$ \\ \hline
\end{tabular}
\normalsize

The corresponding normal operators in the unperturbed situation is as 
below:

\footnotesize

\begin{tabular}{|c|c|c|}\hline
&    Reducible& Normal Operator \\ 
\hline
  SW & Trivial &  $N^{\sw,0}_{\mTheta}(\gamma,a,\phi)=(d^{*}a, *da+d\gamma, 
  (D\oplus D)\phi)$ \\ \cline{2-3}
    &Type I & $N^{\sw,I}_{A}(\phi)=D_{A}\phi$\\ \cline{2-3}
    &Type II &  $N^{\sw,II}_{\alpha,\psi}(\gamma,a,\phi)=
    (d^{*}_{\alpha}a-B(\phi,\psi),*d_{\alpha}a-\{\phi\cdot\psi\}_{0*}+
    d_{\alpha}\gamma,$\\ 
&&$ D_{\alpha}\phi+a\cdot\psi-\gamma\psi)$ \\ \hline
$SU(3)$ &Trivial & $N^{\su,0}_{\mTheta}(\gamma,a)=(d^{*}a,*da+d\gamma)$\\ 
\cline{2-3}
&Type I &   $N^{\su,I}_{\alpha}(\gamma,a)=(d^{*}_{\alpha}a,
*d_{\alpha}a+d_{\alpha}\gamma)$ \\ \hline
\end{tabular}
\normalsize

\bigskip

\proofof{Lemmas~\ref{det-orient-lem}, \ref{det-red-lem}}
Assume Proposition~\ref{spec-prop}. 
The lemmas follow from the unperturbed situation $\pi=0$ by application of the deformation
$t\pi$, $0\le t\le 1$. The assertions of Lemmas~\ref{det-orient-lem}
and \ref{det-red-lem} are then equivalent to the condition that the 
respective operators (which are now formally self-adjoint elliptic)
have even real spectral flow around any path that is a loop at the level of 
the quotient space. Let us now deal with the individual operators in 
turn. 

$K^{\sw,I}_{A}$ and $K^{\su,I}_{A}$ are the same operator, 
the (negative of the) boundary of the Anti-Self-Dual (ASD) operator in 4-dimensions. By 
\cite{APS} the spectral-flow around closed loops is equal to the 
(negative of the) index of a twisted ASD operator
on $Y\times S^{1}$. This index is well-known to be congruent to 
$0\bmod 8$. After a deformation we can decompose $K^{\sw,II}_{\alpha,\psi}$ into a sum 
of three operators: $K\oplus D_{\alpha}\oplus K$ where 
$K(\gamma,a)=(d^{*}a,*da+d\gamma)$. The Dirac operator $D_{\alpha}$ is 
clearly complex and thus even spectral flow. Being topological, $K$ 
has no spectral-flow. Thus $K^{\sw,II}_{\alpha,\psi}$ has even real spectral flow 
around loops in the quotient space. A deformation of $L^{\sw}_{A,\mPhi}$ 
brings it into a direct sum $K^{\sw,I}_{A}\oplus N^{\sw,I}_{A}$. The
complex linear nature of $N^{\sw,I}_{A}$ means it always has even 
real spectral flow whereas the spectral-flow of $K^{\sw,I}_{A}$ we have
already treated. The case of $L^{\su}_{A}$ follows from the same line of 
reasoning as above, being the negative of the boundary of the ASD 
operator. \qed

\proofof{Lemma~\ref{sf-lem}} Assuming Proposition~\ref{spec-prop} 
and after a deformation we may again assume $\pi=0$. Item(a): The path 
$\hat{\gamma}$ defines a $U(2)$-connection $\widehat{A}$ over $Y\times[0,1]$.  
Since $N^{\su,I}_{A}$ is the \textit{negative} of the boundary of the 
$\widehat{A}$-twisted ASD 
operator on $Y\times[0,1]$ (with the orientation $dydt$), by \cite{APS} 
the spectral flow of $N^{\su,I}_{A}$ along $\hat{\gamma}$ is equal to 
the \textit{negative} of the index of the ASD-operator on $Y\times[0,1]$ with 
APS spectral boundary conditions. This index is computed from 
\cite{AHS}, \cite{APS} to be
\begin{equation}
    -2\int_{Y\times[0,1]}c_{2}(\widehat{A}) 
    =-2\Bigl(\cs(\hat{\gamma}(0))-\cs(\hat{\gamma}(1))\Bigr).
    \end{equation}
(Note: with the orientation $dydt$ the boundary of $Y\times[0,1]$ is 
oriented according to $Y\times\{0\}-Y\times\{1\}$ for Stokes' Theorem 
to hold without any signs.)
On the other hand $N^{\sw,I}_{A}$ is the boundary of the 4-dimensional 
Dirac operator coupled to $\widehat{A}$ on $Y\times[0,1]$. Thus the 
spectral flow along $\hat{\gamma}$ is equal to the index of this 
4-dimensional Dirac operator.
This has 
index given by $-\int_{Y\times[0,1]}c_{2}(\widehat{A})$.
Item(b) can be seen by either a similar computation or the 
observation that $\CC^{0,II}/\G_{E}$ is simply connected. \qed

\proofof{Proposition~\ref{spec-prop}} (\cite{taubes}) we need to show that 
$N^{\pi}_{x}$ is essentially self-adjoint and has compact resolvent; 
then the spectrum is real and discrete. The absence of
accumulation points and
unboundedness as $\to\pm\infty$ follows by the Hilbert-Schmidt 
theorem applied to 
any resolvent of $N^{\pi}_{x}$. In the 
unperturbed case, $N_{x}$ is formally self-adjoint elliptic on a 
compact manifold and standard elliptic theory gives $N_{x}$ as 
essentially self-adjoint and with compact resolvent. In general we can
regard $N^{\pi}_{x}$ as a perturbation of $N_{x}$ since
$N^{\pi}_{x}=N_{x}+(L^{\nu}\pi)_{x}$. By Lemma~\ref{nor-sym-prop} 
$(L^{\nu}\pi)_{x}$ is symmetric 
with dense domain the $L^{2}_{2}$-sections, and thus is closable. 
Furthermore since $(L^{\nu}\pi)_{x}$ is bounded as an operator on 
$L^{2}_{2}$-sections it follows that $(L^{\nu}\pi)_{x}$ is relatively 
compact with respect to $N_{x}$. This  in turn implies that $(L^{\nu}\pi)_{x}$
has arbitarily small relative bound with respect to $N_{x}$. Now standard 
stability theory \cite{kato} tells us $N^{\pi}_{x}$ is 
also essentially self-adjoint and has compact resolvent. \qed

\subsection{The Main Theorem}

Recall $E\to Y$ is our $U(2)$-bundle in SW-theory. Without loss, let us choose
now in $SU(3)$-Casson $F=E\oplus \overline{\det E}$. Then given a connection
$A$ on $E$, this induces the $SU(3)$ connection $A\oplus\overline{\det A}$ on $F$.
In this way we obtain an identification
\begin{equation}\label{red-identify}
\CC^{0,I}_{E}/\G_{E} = \A^{0,I}_{F}/\G_{F}.
\end{equation}
Henceforth we shall assume this identification.

\begin{lemma}
There exists admissible non-degenerate perturbations $\pi$, $\pi'$ 
in $U(2)$-SW, $SU(3)$-Casson-Taubes respectively such that
$\M^{\sw,I}_{\pi}\cong\M^{\su,I}_{\pi'}$. 
\end{lemma}

\proof Follow proof of Proposition~\ref{adm-lem-sw}. \qed

Such a pair of perturbations we call \textit{non-degenerate 
compatible}. In such a situation we shall simply write 
$\M^{I}_{\pi}$ for both $\M^{\sw,I}_{\pi}$
and $\M^{\su,I}_{\pi'}$, regarding them as being 
identified.

At a non-degenerate point $x\in \M^{\sw*}_{\pi}$ ($\M^{su*}_{\pi'}$) the
kernel and cokernel of $L^{\sw,\pi}_{\overline{x}}$ ($L^{\su,\pi'}_{\overline{x}}$)
are trivial, $\overline{x}$ any representative of $x$. 
This gives a canonical trivialization (with $\R$) of the determinant
line $\detind \widehat{L}^{\sw,\pi}$ ($\detind \widehat{L}^{\su,\pi'}$) at $x$.
Thus in order to orient $\M^{\sw*}_{\pi}$ ($\M^{\su*}_{\pi'}$) we
fix the overall orientation of $\detind \widehat{L}^{\sw,\pi}$ ($\detind 
\widehat{L}^{su,\pi'}$)
by fixing the orientation at one point, which we take to be $[\mTheta]$
and then propagating the orientation from this point. This is well-defined
by Lemma~\ref{det-orient-lem}.
At $\overline{x}=\mTheta$ the kernel and cokernel of $L^{\sw,\pi}_{\overline{x}}$ or 
$L^{su,\pi'}_{\overline{x}}$
are identical, call them $H$ and specify the orientation of the 
determinant line at $[\mTheta]$ by the
rule $o(H)\wedge o(H)^{*}$ where $o(H)$ is any orientation of $H$ and 
$o(H)^{*}$ the dual orientation. We denote the orientation at $x$ by
$$
\epsilon(x) \in \{\pm 1\}.
$$
In an identical manner non-degenerate points in the reducible strata 
$\M^{\sw,I}_{\pi}$, $\M^{\sw,II}_{\pi}$, $\M^{\su,I}_{\pi'}$
are oriented but this time using the tangential
operators $K^{\sw,I,\pi}_{A}$, 
$K^{\sw,II,\pi}_{A_{0},\mPhi_{0}}$,
$K^{\su,I}_{A}$ respectively and invoking Lemma~\ref{det-red-lem}.
Again we denote the orientation at $x$ by $\epsilon(x)$, the usage will be
clear from the context.

In the situation of compatible perturbations,
the operators $K^{\sw,I,\pi}_{A}$
and $K^{\su,I,\pi'}_{A}$ coincide. Thus we identify them
and simply write $K^{I,\pi}_{A}$. It is clear that the 
orientation at a non-degenerate point $x\in\M^{I}_{\pi}$
is the same in the $U(2)$-SW and $SU(3)$-Casson theories.

Assume  henceforth that $\pi$, $\pi'$ are non-degenerate compatible
perturbations. For the top (non-singular) strata contribution, let
\begin{equation}\label{irred-term}
\Lambda^{*}_{\sw}(g,\pi) = \sum_{x\in\M^{\sw*}_{g,\pi}}\epsilon(x),\quad
\Lambda^{*}_{\su}(g,\pi') =\sum_{x\in\M^{\su*}_{\pi'}}\epsilon(x).
\end{equation}

\textbf{Counter-terms (0)}
We remind the reader that $g$ denotes the metric on our ZHS $Y$ and
$D$ the canonical Dirac operator on the
spinor bundle $S\to Y$. Denote by $B$ the operator on $Y$ which is 
half the boundary of the signature operator in 4-dimensions. To these operators we may
associate the Atiyah-Patodi-Singer (APS) spectral invariants \cite{APS}
\begin{equation}
\eta(B),\quad \xi = \frac{1}{2}\Bigl( \eta(D) + \dimn_{C}\ker D\Bigr).
\end{equation}
Let $\pi$ be a perturbation in $U(2)$-SW-theory and let $D^{\pi}_{\mTheta}$ the
normal operator at a trivial connection $\mTheta$. Recall that this operator
acts on $S\otimes E$.
Now set 
\begin{eqnarray}
c(g,\pi) &=& \eta(B) +\frac{1}{8}\xi 
+\frac{1}{2}\Bigl(\mbox{complex spectral flow of}\\ 
&&\phantom{\eta(B) +\frac{1}{8}\xi+\frac{1}{2}\Bigl(}\mbox{$\{ (1-t)D_{\mTheta}
 +tD^{\pi}_{\mTheta} \}^{1}_{t=0}$}\Bigr)\nonumber
\end{eqnarray}

\begin{lemma}
$c(g,\pi) \equiv \mu(Y) \bmod 2$
where $\mu(Y)\in \{0,1\}$ is the Rokhlin invariant.
\end{lemma} 

\proof This is discussed in \cite{lim0} but we repeat it here
for convenience. 
If $X$ is compact oriented spin 4-manifold with oriented boundary $Y$
then an application of the APS index theorems to $X$ shows that
\begin{equation}\label{ind-equ}
\xi+\frac{1}{8}\eta(B) = -\mbox{Index}\, D^{(4)} -\frac{1}{8}\mbox{sign}\,X.
\end{equation}
Here $D^{(4)}$ is the Dirac operator on $X$ and ${\rm sign}\,X$ the signature.
Thus we see that the left-side of \ref{ind-equ} is always an integer.
The $\bmod 2$ reduction of the right-side only involves the
signature term (since in four dimensions the Dirac operator is quaternionic linear and
so its index is even) and therefore is just the Rokhlin invariant $\mu(Y)$. 
The complex spectral flow of the family $\{ (1-t)D_{\mTheta}
 +tD^{\pi}_{\mTheta}\}^{t=1}_{t=0}$
is always divisible by 4 since both $D_{\mTheta}$ and $D^{\pi}_{\mTheta}$ are 
each a double of a quaternionic linear operator. \qed

Now set, as the contribution from the SW-strata of the trivial connection:
\begin{equation}
\Lambda^{0}(g,\pi) = \frac{1}{2}c(g,\pi)\Bigl( c(g,\pi) -1\Bigr).
\end{equation}

\textbf{Counter-terms (I)}
For the Type I reducible strata contribution set
\begin{equation}\label{red-I-term}
\Lambda^{I}(g,\pi,\pi') = \sum_{x\in\M^{I}_{\pi}}
\epsilon(x)
\Bigl\{(\SF^{\su,I}_{\nu}+2\SF^{\sw,I}_{\nu})([\mTheta],x)+4c(g,\pi')\Bigr\}
\end{equation}

The spectral flow term is taken along any path in $\CC^{0,1}_{E}/\G_{E}=
\A^{0,I}_{F}/\G_{F}$ from $[\mTheta]$ to $x$. This 
is  well-defined independent of path by Lemma~\ref{sf-lem}.

\textbf{Counter-terms (II)} Lastly, for the SW-Type II reducible strata set
\begin{equation}\label{red-II-term}
\Lambda^{II}(g,\pi) = \sum_{x\in\M^{\sw,II}_{\pi}}
\epsilon(x)\Bigl\{ \SF^{\sw,II}_{\nu}([\mTheta],x) + c(g,\pi)\Bigr\}.
\end{equation}
This is again well-defined, by Lemma~\ref{sf-lem}.

As a remark, the presence of the metric dependent terms $c(g,\pi)$ in
$\Lambda^{I}$ and $\Lambda^{II}$ are inserted to cancel out spectral 
flow phenomena at $[\mTheta]$ in expressions $\SF^{\sw,I}_{\nu}([\mTheta],x)$,
$\SF^{\sw,II}_{\nu}([\mTheta],x)$ if we were to vary the metric.

The main result of this paper is:

\begin{theorem}\label{main-th}
Let $Y$ be an oriented closed integral homology 3-sphere with Riemannian metric $g$.
Let $\pi$, $\pi'$ be non-degenerate compatible admissible perturbations in $U(2)$-SW,
$SU(3)$-Casson respectively. Then the sum 
\begin{displaymath}
\tau(Y) = \Lambda^{*}_{\su}(g,\pi')+2\Lambda^{*}_{\sw}(g,\pi)
+\Lambda^{I}(g,\pi,\pi')+2\Lambda^{II}(g,\pi)+2\Lambda^{0}(g,\pi)
\end{displaymath}
is an integer and
independent of $g$ and $\pi$, $\pi'$ and thus defines an oriented 
diffeomorphism invariant of $Y$.
\end{theorem}

Let
\begin{displaymath}
    \lambda^{\su(2)}_{g,\pi}(Y)=\sum_{x\in\M^{I}_{\pi}}\epsilon(x),
    \quad
    \lambda^{\sw}_{g,\pi}(Y)=\sum_{x\in\M^{\sw,II}_{\pi}}\epsilon(x) +c(g,\pi).
    \end{displaymath}
According to Taubes \cite{taubes}, $\lambda^{\su(2)}_{g,\pi}(Y)$ is 
(up to a universal sign) twice Casson's invariant. 
$\lambda^{\sw}_{g,\pi}(Y)$ on the other hand is precisely the 
definition of the abelian SW-invariant for $Y$ (see for instance 
\cite{lim0}). By \cite{lim1}, this is again (up to a universal 
sign) equal to Casson's invariant.

\begin{cor}\label{reverse}
    Let $-Y$ denote $Y$ but with the reverse orientation. Then
    \begin{displaymath}
	 \tau(-Y)=\tau(Y)+2\lambda^{\su(2)}_{g,\pi}(Y)+2\lambda^{\sw}_{g,\pi}(Y).
	 \end{displaymath}
Thus an orientation independent invariant for $Y$ is given by
\begin{displaymath}
    \overline{\tau}(Y) 
    =\tau(Y)+\lambda^{\su(2)}_{g,\pi}(Y)+\lambda^{\sw}_{g,\pi}(Y).
   \end{displaymath}
\end{cor}

\proof
Under reversal of orientation of $Y$ the $\pi'$-perturbed flat equation 
transforms to the $(-\pi')$-perturbed equations. Thus we can identify 
the $SU(3)$-moduli spaces in either orientation.
The Floer-Taubes operator 
switches to its negative. In
the SW-case, reversal simply changes the Clifford action on the 
spinor bundle to its negative. Thus the SW-equation in the reversed 
orientation is equation to the original except that the Dirac 
component of the equation switches to its negative. If $\pi=(*k,l)$ 
is the non-degenerate perturbation originally used then 
$\overline{\pi}=(*k,-l)$ is non-degenerate and admissible in the 
reversed context. If $(A,\mPhi)$ is a (perturbed) solution in 
the original then $(A,-\mPhi)$ is a (perturbed) solution in the 
reversed. Thus we can identify $\M^{\sw}_{\overline{\pi}}(-Y)$ with $\M^{\sw}_{\pi}(Y)$.
The Floer-Taubes operator also changes to its negative in the 
SW-context.

Under a reversal of orientation $Y\mapsto -Y$ we obtain the 
transformations
\begin{eqnarray*}
    \epsilon(x) &\mapsto& \epsilon(x), 
    x\in\M^{\sw*}_{\pi},\M^{\su*}_{\pi'}\\
    \epsilon(x) &\mapsto& -\epsilon(x), 
    x\in\M^{\sw,r}_{\pi},\M^{\su,r}_{\pi'}\\
    \SF_{\nu}^{\su,I}([\mTheta],x) &\mapsto& -\SF_{\nu}^{\su,I}([\mTheta],x)-2\\
    \SF^{\sw,I}_{\nu}([\mTheta],x) &\mapsto& 
    -\SF^{\sw,I}_{\nu}([\mTheta],x)\\
    \SF^{\sw,II}_{\nu}([\mTheta],x) &\mapsto& 
    -\SF^{\sw,II}_{\nu}([\mTheta],x)-1\\
    c(g,\pi) &\mapsto& -c(g,\pi).
    \end{eqnarray*}
The non-trivial constants in the spectral flow terms are the 
corrections terms which equal the dimension of the kernel of the normal 
operators at $\mTheta$.
The orientation reversal formula easily follows.  
$\overline{\tau}(Y)$ is the average of $\tau(Y)$ and $\tau(-Y)$ and 
thus independent of orientation. 
\qed

\begin{remark}\rm
    If the Casson invariant of $Y$ is non-zero then by 
    Corollary~\ref{reverse} at least one of 
$\tau(Y)$ or $\tau(-Y)$ is also non-zero. Since there are infinitely 
many $Y$ for which Casson's invariant is non-zero we deduce that
there are infinitely many integral homology spheres for which
$\tau\neq 0$  for at least one orientation.
\end{remark}


\section{Proof of Theorem~\protect\ref{main-th}}\label{proof-main}

The basic strategy of the proof is an extension of those in \cite{lim0},
\cite{BH} and \cite{lim3}. Some (standard) portions of the argument 
are omitted and can be found in the cited references.

\textbf{Standing Convention} To simplify notation we often confuse a point say 
$x$, in the quotient space $\ZZ/\G$ with a representative of $x$ in 
the context of the Floer-Taubes/normal/tangential operators as 
well as cohomology spaces. This is permissible as different 
representatives of the same gauge orbit give rise to conjugate 
operators. In particular we may regard the normal operators as being 
parameterized by $\ZZ^{r}/\G$.

\subsection{Parameterized Moduli Spaces}

Let $(g_{0},\pi_{0},\pi_{0}')$ and $(g_{1},\pi_{1},\pi_{1}')$ be triples
consisting of metric and non-degenerate admissible compatible
perturbations. We wish to compare the moduli spaces for these two triples.

To this end we may form the corresponding \textit{parameterized 
moduli spaces} in the SW and SU(3) cases:
\begin{eqnarray}
W^{\sw}=\bigcup_{t\in[0,1]}\M^{\sw}_{g_{t},\pi_{t}}\times\{t\}\subset\CC_{E}/\G_{E}\times[0,1],\\
W^{\su}=\bigcup_{t\in[0,1]}\M^{\su}_{g_{t},\pi_{t}}\times\{t\}\subset\A_{F}/\G_{F}\times[0,1].\nonumber
\end{eqnarray}
(Note: in the SW-case, fix a model for the spinor bundle. Then regard the 
Clifford action etc., as varying with the metric.)
We retain the usage of the
(superscript) notations $*,r,I,II$ etc.  pertaining to the unparameterized moduli spaces  
in the parameterized setting. The parameterized moduli spaces can also be 
regarded as being formed from the $t$-dependent SW and SU(3)-flat 
equations over $\ZZ\times[0,1]$:
\begin{equation}\label{para-chi}
    \widetilde{\XX}^{\sw}(A,\mPhi,t)=0, \quad 
    \widetilde{\XX}^{\su}(A)=0.
    \end{equation}
It will be necessary to introduce perturbations in the parameterized 
context.
An admissible perturbation $\sigma$ in the 
parameterized context is defined in an analogous manner to the 
unparameterized case as a 
function $\sigma\colon\ZZ\times [0,1]\to L^{2}_{2}(V)$
with the additional condition that the support is
contained in $\ZZ\times(0,1)$. The restriction of $\sigma$ to a slice
$\ZZ\times\{t\}$ is clearly an admissible perturbation on $Y$ itself.

Denote the $\sigma$, $\sigma'$ perturbed versions of the parameterized moduli spaces by 
\begin{equation}
W^{\sw}_{\sigma},\quad W^{\su}_{\sigma'}
\end{equation}
respectively. To save notation we shall use $W_{\sigma}$ to denote either 
space.

\begin{prop} 
$W_{\sigma}$  is always compact.
\end{prop}

\proof Follows from Proposition~\ref{compact-prop}. \qed

In compact notation we may write the $\sigma$-perturbed (\ref{para-chi}) as
\begin{equation}\label{para-map1}
    \widetilde{\XX}_{\sigma}\colon\ZZ\times[0,1]\to L^{2}_{2}(V).
\end{equation}
If we let the linearization be $\widetilde{\delta}^{1}_{x,t}$ then 
this fits into an elliptic complex
\begin{equation}\label{para-ellip}
    L^{2}_{3}(W)\stackrel{\delta^{0}_{x}}{\longrightarrow}
    L^{2}_{2}(V)\oplus\R\stackrel{\widetilde{\delta}^{1}_{x,t}}{\longrightarrow}
    L^{2}_{1}(V)\stackrel{\delta^{0*}_{x}}{\longrightarrow}
    L^{2}(W)
\end{equation}
with (harmonic) cohomologies $\widetilde{\bH}^{i}_{x,t}$.
The condition that $W^{*}_{\sigma}$ is cut out
transversely is $\widetilde{\bH}^{1}_{x,t}=\{0\}$, i.e. 
non-degenerate/regular.   

Along a reducible strata again we have an orthogonal decomposition 
paralleling Proposition~\ref{ellip-split}, with tangential and normal 
cohomologies $\widetilde{\bH}^{i,\tau}_{x,t}$ and 
$\widetilde{\bH}^{i,\nu}_{x,t}$ respectively. Call 
$W^{\alpha}_{\sigma}$, $\alpha\in\{0,I,II\}$ \textit{regular} if
$\widetilde{\bH}^{1,\tau}_{x,t}=\{0\}$, $(x,t)\in W^{\alpha}_{\sigma}$. 
This is the condition of being cut out \textit{transversely in 
$\ZZ^{\alpha}/\G\times[0,1]$}.

In the next proposition let $\M$ denote either $\M^{\sw}$ 
or $\M^{\su}$.

\begin{prop}\label{struct-prop}
There are admissible perturbations $\sigma$  such that 
$W^{\sw}_{\sigma}$ is a regular stratified compact
singular cobordisms between $\M_{\pi_{0}}$ and $\M_{\pi_{1}}$, i.e. the following hold:
\newline
{\rm(a)} Each individual strata $W^{\alpha}_{\sigma}$, 
$\alpha\in\{*,0,I,II\}$ is regular: it is a 1-manifold with boundary $\M_{\pi_{0}}\cup\M_{\pi_{1}}$ and with
possibly a number of non-compact ends. 
\newline
{\rm(b)} Each end limits to a {\rm singular point}, which 
lies in a  reducible strata. There are only finitely many
singular points. 
\newline
{\rm(c)}  A neighbourhood of each singular point is diffeomorphic to 
$T=\{(r,t)\,|\, rt=0, t\ge 0\}\subset\R^{2}$ where the edge 
$\{0\}\times(0,\infty)$ corresponds to the limiting end.
\newline
{\rm(d)} Each reducible strata parameterizes the associated
family of normal operators; the singular points are exactly where the
family experiences spectral flow. 
The spectral flow at these points are always transverse
\newline
{\rm(e)} Limiting ends only occur in the irreducible and in the SW-Type II strata.
When a singular point lies in the SW-Trivial strata the corresponding
limiting end lies in the SW-Type II reducible strata, and every 
limiting end of the SW-Type II strata lies in the SW-Trivial strata.
There are no singular points on the SU(3)-Trivial strata.
\end{prop}

This will be proven in \S\ref{sub-proof}.

\begin{remark} \rm
    (a) The normal operator at $(x,t)\in W^{r}_{\sigma}$ is the normal operator at 
    $x$ in $\M^{r}_{g_{t},\pi_{t}}$
    (b) Transverse spectral flow
means all eigenvalues crosses zero transversely, modulo
multiplicity if the operator has a complex or quaternionic structure, i.e.
we have only  simple eigenvalues over $\R$, $\C$ or $\qu$.
\end{remark}

\begin{prop}\label{orient-para-lem}
Let $\sigma,\sigma'$ be perturbations as in 
Proposition~\ref{struct-prop}.
$W^{\sw}_{\sigma}$ and $W^{\sw}_{\sigma'}$ admit a consistent orientation
convention such that 
\begin{eqnarray*}
\partial (W^{\sw}_{\sigma}\backslash\{\mbox{\rm singular points}\}) 
&=& \M^{\sw}_{\pi_{1}}\cup -\M^{\sw}_{\pi_{0}} \\
\partial (W^{\su}_{\sigma'}\backslash\{\mbox{\rm singular points}\}) 
&=& \M^{\su}_{\pi_{1}'}\cup -\M^{\su}_{\pi_{0}'}.
\end{eqnarray*}
\end{prop}
Note that the assertion of the proposition is inclusive of the reducible stratas
where we assign the orientation value $+1$ to the trivial connection 
$[\mTheta]$.
The existence of the claimed orientations is established by
considering determinant line bundles and will be established
in the section on orientation, \S\ref{sec-orient}.

The existence of the orientations on the reducible strata of the
parameterized moduli spaces also allows us to assign to each singular
point a value $+1$ or $-1$ according to the spectral flow of the normal
operator at that point, moving in the direction of the given orientation.

\begin{prop}\label{sing-pt-lem}
Let $\mXi^{-1}(0)$, $\mXi\colon \R\times[0,\infty)\to\R$, $\mXi(x,y) = xy$
 be a orientation preserving local model for a singular point, where the orientation on
 $\R\times\{0\}$ is the usual orientation on $\R$. Let $\epsilon\in\{\pm 1\}$
 be the sign of the spectral flow of the normal operator at $x=0$ in the local
 model. Then the orientation on $\{0\}\times [0,\infty)$ is given by
 $-\epsilon$ multiplied with the standard orientation on $[0,\infty)$.
 \end{prop}

The proof is also in \S\ref{sec-orient}.

\subsection{The Main Argument}
 
To show that the sum $\tau(Y)$ in Theorem~\ref{main-th} is a diffeomorphism
invariant we need to show that the defect
\begin{eqnarray}\label{irred-defect}
2\Lambda^{*}_{\sw}(g_{1},\pi_{1},\pi_{1}')
+\Lambda^{*}_{\su}(g_{1},\pi_{1},\pi_{1}')
-2\Lambda^{*}_{\sw}(g_{0},\pi_{0},\pi_{0}')
-\Lambda^{*}_{\su}(g_{0},\pi_{0},\pi_{0}')
\end{eqnarray}
exactly cancels the defects
\begin{eqnarray}
&\Lambda^{I}(g_{1},\pi_{1},\pi_{1}')-\Lambda^{I}(g_{0},\pi_{0},\pi_{0}'),\nonumber\\
&2\Lambda^{II}(g_{1},\pi_{1})-2\Lambda^{II}(g_{0},\pi_{0})\label{red-defect},\\
&2\Lambda^{0}(g_{1},\pi_{1})-2\Lambda^{0}(g_{0},\pi_{0}).\nonumber
\end{eqnarray}
This is established by through the singular cobordisms 
$W^{\sw}_{\sigma}$ and $W^{\su}_{\sigma'}$ of Proposition~\ref{struct-prop}.
To this end, without loss we may assume that
$W^{\sw}_{\sigma}$ and $W^{\su}_{\sigma'}$ are \textit{elementary}
singular cobordisms, by which we mean the occurance of exactly one
singular point (or none). 
According to
Proposition~\ref{struct-prop} we have the following different types of
elementary singular cobordisms for $W^{\sw}_{\sigma}$ and $W^{\su}_{\sigma'}$.
\begin{itemize}
\item no singular points
\item singular point is Type I reducible
\item singular point is Type II reducible
\item singular point is Trivial (connection)
\end{itemize}
In the case of all these types of elementary singular cobordisms
except for the last (let us call it $W^{\sw,0}$),
the invariance by analysing the defects, is  covered in \cite{BH}
(see also \cite{lim3}, \cite{lim0}) without any new idea. These cases
correspond to the \textit{birth} or \textit{death} of new points at bifurcation
into the highest (i.e. irreducible) strata in the parameterized moduli space.
The defect (\ref{irred-defect}) is cancelled by the first two defects in
(\ref{red-defect}) with the last defect in (\ref{red-defect}) 
identically zero.

The case $W^{\sw,0}$ presents the new phenomena of birth/death of
new points into the Type II strata at bifurcation. As such it represents
a \lq second order\rq\ defect, being the
defect of the counter-term $\Lambda^{I}(g,\pi,\pi')$. 
To simply matters more, we may assume that $W^{\sw,0}$ consists of components all of which
are topologically closed intervals $[0,1]$ except a single one which is
topologically
$[-1,1]\times\{0\}\cup\{0\}\times[0,1]$. (There are actually two subcases 
corresponding to where the boundary of the normal edge lies.) 
Clearly in this case the  defect
(\ref{irred-defect}) is zero as well as the term
$\Lambda^{I}(g_{1},\pi_{1},\pi_{1}')-\Lambda^{I}(g_{0},\pi_{0},\pi_{0}')$
since the irreducible and Type I strata components are assumed to be
closed intervals $[0,1]$. Thus we only have to deal with the
changes in the defect terms $\Lambda^{II}$ and $\Lambda^{0}$.

\begin{claim}\label{pf-claim}
$\Lambda^{II}(g_{1},\pi_{1},\pi_{1}')-\Lambda^{II}(g_{0},\pi_{0},\pi_{0}')
+\Lambda^{0}(g_{1},\pi_{1})-\Lambda^{0}(g_{0},\pi_{0})=0$.
\end{claim}

Given the claim, Theorem~\ref{main-th} is proven.

\subsection{Proof of Claim~\protect\ref{pf-claim}}

For convenience we change notation; assume the parameterization varies over 
$[-1,1]$ instead of $[0,1]$ so the initial metric, perturbation, etc. 
are now $g_{-1}$, $\pi_{-1}$ etc.

Any component of $W^{\sw,0}$ which is a product makes no contribution to
the defects $\Lambda^{II}(g_{1},\pi_{1},\pi_{1}')-\Lambda^{II}(g_{-1},\pi_{-1},\pi_{-1}')$
and $\Lambda^{0}(g_{1},\pi_{1})-\Lambda^{0}(g_{-1},\pi_{-1})$ so we 
focus our attention on the component of $W^{\sw,0}$ which topologically
is $[-1,1]\times\{0\}\cup\{0\}\times[0,1]$. Let $\epsilon=\pm 1$ be the sign
of the $\qu$-spectral flow (or half the $\C$-spectral flow) of the normal operator 
$N^{\sw,0}$at the singular point $(0,0)$. By Lemma~\ref{sing-pt-lem}
the arc $\{0\}\times [0,1]$ is oriented as $-\epsilon$ times the standard orientation
on $[0,1]$. We have two situations for the boundary point $p=(0,1)$. Denote
by Case A when this is in $\M^{\sw}_{\pi_{-1}}$  and
Case B when in $\M^{\sw}_{\pi_{1}}$. Recall that $\epsilon(p)\in\{\pm 
1\}$ denotes the orientation of $p$ as a point in 
$\M^{\sw}_{\pi_{-1}}$ or $\M^{\sw}_{\pi_{1}}$.
Observe in Case A, 
$\epsilon(p)=\epsilon$ and in Case B, $\epsilon(p)=-\epsilon$.
Then the defect
\begin{eqnarray*}
\lefteqn{\Lambda^{II}(g_{1},\pi_{1},\pi_{1}')-\Lambda^{II}(g_{-1},\pi_{-1},\pi_{-1}')}\\
&=&\left\{
\begin{array}{ll}
-\epsilon\Bigl( \SF^{\sw,II}_{\nu}([\mTheta],p) +c(g_{-1},\pi_{-1})\Bigr)&
\mbox{in Case A}\\
-\epsilon\Bigl( \SF^{\sw,II}_{\nu}([\mTheta],p) +c(g_{1},\pi_{1})\Bigr)&
\mbox{in Case B}.
\end{array}\right.
\end{eqnarray*}
In the notation we implicitly assume that $\SF^{\sw,II}_{\nu}([\mTheta],p)$ takes 
place in either $\CC^{r}_{E}\times\{-1\}$ or $\CC^{r}_{E}\times\{+1\}$ 
depending on where $p$ is located.
The key observation is the following:

\begin{lemma}\label{claim-lem}
\begin{displaymath}
\SF^{\sw,II}_{\nu}([\mTheta],p) =\left\{
\begin{array}{ll}
\frac{1}{2}(\epsilon-1)&\mbox{\rm in Case A}\\
-\frac{1}{2}(\epsilon+1)&\mbox{\rm in Case B}
\end{array}
\right.
\end{displaymath}
\end{lemma}

This shall be proven below. It follows then that
\begin{eqnarray*}
\lefteqn{\Lambda^{II}(g_{1},\pi_{1},\pi_{1}')-\Lambda^{II}(g_{-1},\pi_{-1},\pi_{-1}')}\\
&=&\left\{
\begin{array}{ll}
-\epsilon\Bigl( \frac{1}{2}(\epsilon-1) +c(g_{-1},\pi_{-1})\Bigr)&
\mbox{in Case A}\\
-\epsilon\Bigl( -\frac{1}{2}(\epsilon+1) +c(g_{1},\pi_{1})\Bigr)&
\mbox{in Case B}
\end{array}\right.\\
&=& \frac{1}{2}(\epsilon-1)-\epsilon c(g_{-1},\pi_{-1}),
\end{eqnarray*}
where in the last line we use the relation $c(g_{1},\pi_{1})=c(g_{-1},\pi_{-1})+\epsilon$.
To prove this recall from the definition that $c(g_{t},\pi_{t})$ changes by the spectral flow 
of $D^{\pi_{t}}$ (with respect to metric $g_{t}$) acting on $S$, as $t$ 
varies. We claim this is 
exactly half the $\C$-spectral flow of $N^{\sw,0}_{\Theta,t}$
as $t$ varies. After trivializing $E$ as $\C^{2}\times Y$ using 
$\mTheta$, it is seen that $N^{\sw,0}_{\Theta,t}=K\oplus 
D^{\pi_{t}}\oplus D^{\pi_{t}}$ where $K$ is the deRham 
operator on $\mLambda^{0+1}\otimes\C$ (see \S\ref{detailproof}). Since $K$ is topological it 
has no spectral flow so half the $\C$-spectral flow of 
$N^{\sw,0}_{\mTheta,t}$ is equal to the $\C$-spectral flow of 
$D^{\pi_{t}}$, as claimed.

To continue: on the other hand we easily see the defect
\begin{eqnarray*}
\lefteqn{\Lambda^{0}(g_{1},\pi_{1})-\Lambda^{0}(g_{-1},\pi_{-1})}\\
&=& \frac{1}{2}\Bigl( c(g_{-1},\pi_{-1})+\epsilon\Bigr) 
\Bigl(c(g_{-1},\pi_{-1})+\epsilon-1\Bigr)\\
&&\mbox{}-\frac{1}{2} c(g_{-1},\pi_{-1})\Bigl(c(g_{-1},\pi_{-1})-1\Bigr)\\
&=& \frac{1}{2}(1-\epsilon)+\epsilon c(g_{-1},\pi_{-1}).
\end{eqnarray*}
Hence the sum of the defects in zero and Claim~\ref{pf-claim}
is established. \qed

\proofof{Lemma~\ref{claim-lem}}
Let $A$ reduce as $\theta\oplus\theta$ in the splitting $E=L_{0}\oplus L_{1}$
and let $\mPhi=(\phi,0)\in L^{2}_{2}(S\otimes L_{0})\oplus 
L^{2}_{2}(S\otimes L_{1})$. Furthermore trivialize $L_{i}$ as $\C\times 
Y$ via the trivial connection $\theta$.
Then the normal operator $N^{\sw,II,\pi}_{\theta,\phi,t}$ acts
on sections of $(\mLambda^{0+1}\otimes\C)\oplus S$ 
(\S\ref{detailproof}).
If $\phi=0$, then 
$N^{\sw,II,\pi}_{\theta,\phi,t}$ decouples as 
$K\oplus D^{\pi_{t}}$ where $K$ is the
deRham operator on $\mLambda^{0+1}\otimes \C$ and
$D^{\pi_{t}}$ is the $\pi_{t}$-perturbed Dirac operator on $S$ with respect to 
metric $g_{t}$.  Thus the
spectral flow of $N^{\sw,II,\pi}$ along the arc $[-1,0]\times\{0\}$
in the local model for the singular point is $\frac{1}{2}(\epsilon+1)$
and along the arc $[0,1]\times\{0\}$ it is $\frac{1}{2}(\epsilon-1)$ (in our
convention, ff. Prop.~\ref{spec-prop}).
Note that the kernel and cokernel of $K$ is $\cong\C$, the constant
functions, and therefore makes no
contribution to spectral flow.

By assumption of transverse spectral flow, $D^{\pi,t}$ has kernel $\cong\C$
at the singular point $(0,0)$. Let $\phi$ be an element, say of unit length
in the kernel. Then to first order, the family 
$N_{s}\colon=N^{\sw,II,\pi}_{\theta,s\phi,0}$,
$s\in[0,1]$ models the family $N^{\sw,II,\pi}$ along $\{0\}\times[0,1]$
at $(0,0)$.

\begin{sublemma}\label{deriv-sub}
Let $N'_{0}$ be the derivative of $N_{s}$ at $s=0$. Identify
$\kernel N_{0}$ with $\C^{2}$ via the basis $\{1,\phi\}$. 
Then the restriction of $N'_{0}$ to $\kernel N_{0}$ followed by $L^{2}$-projection
onto the same is given by the matrix
\begin{equation}
\left(
\begin{array}{cc}
 0&-1\\ -1&0
\end{array}
\right)
\end{equation}
\end{sublemma}

It easily follows from the sublemma that to 1st order the complex eigenvalues
of the family $N_{s}$ at $s=0$ and therefore $N^{\sw,II,\pi}$ along $\{0\}\times[0,1]$
at $(0,0)$ is given by 
\begin{displaymath}
\left|
\begin{array}{cc}
-\lambda&-s\\-s&-\lambda
\end{array}
\right| =0,
\end{displaymath}
i.e. $\lambda=+s$ and $\lambda=-s$. Since along $\{0\}\times[0,1]$
there is no other spectral flow (by assumption), the spectral flow must actually be $-1$. Thus
in Case A,
\begin{displaymath}
\SF^{\sw,II}_{\nu}([\mTheta],p) = \frac{1}{2}(\epsilon+1) -1 = \frac{1}{2}(\epsilon-1)
\end{displaymath}
and in Case B,
\begin{displaymath}
\SF^{\sw,II}_{\nu}([\mTheta],p) = -\frac{1}{2}(\epsilon-1)-1 =-\frac{1}{2}(\epsilon+1).
\end{displaymath}
This proves the Lemma~\ref{claim-lem}, modulo Sublemma~\ref{deriv-sub}. \qed

\proofof{Sublemma~\ref{deriv-sub}}
Our goal is to get an explicit expression for $N_{s}$. Let 
$\delta^{i,\nu}_{\theta,s\phi}$ denote the 
differentials in the normal component of orthogonal decomposition of the fundamental elliptic 
complex along $\CC^{II}(L_{0},L_{1})$ (Proposition~\ref{ellip-split}) 
extended over the trivial strata. Recall $L_{i}$ is trivialized as 
$\C\times Y$ by $\theta$.
Then we have
\begin{eqnarray}\label{delta-0-equ}
    &\delta^{0,\nu}_{\theta,s\phi}\colon 
    L^{2}_{3}(\mLambda^{0}\otimes\C)\to 
    L^{2}_{2}(\mLambda^{1}\otimes\C)\oplus L^{2}_{2}(S) \\
    &\xi \mapsto (d\xi, -\xi (s\phi)),\nonumber
\end{eqnarray}
\begin{eqnarray}\label{delta-1-equ}
    &\delta^{1,\nu}_{\theta,s\phi}\colon
    L^{2}_{2}(\mLambda^{1}\otimes\C)\oplus L^{2}_{2}(S) 
    \to L^{2}_{1}(\mLambda^{1}\otimes\C)\oplus L^{2}_{1}(S)\\
    &(a,\psi)\mapsto (*da-*\{\psi\cdot 
    s\phi\}_{0},D^{\pi,0}\psi+a\cdot(s\phi)).\nonumber
\end{eqnarray}
The $L^{2}$-adjoint of (\ref{delta-0-equ}) is then
\begin{eqnarray}\label{delta-0*-equ}
    &\delta^{0,\nu *}_{\theta,s\phi}\colon 
    L^{2}_{2}(\mLambda^{1}\otimes\C)\oplus L^{2}_{2}(S)\to
    L^{2}_{1}(\mLambda^{0}\otimes\C)\\
    &(a,\psi)\mapsto d^{*}a-\langle{\psi,s\phi}\rangle_{\C}.\nonumber
\end{eqnarray}
$N_{s}$ is given by the block matrix
\begin{displaymath}
    \left(
    \begin{array}{cc}
        0 & \delta^{0,\nu *}_{\theta,s\phi}  \\
	\delta^{0,\nu}_{\theta,s\phi} & \delta^{1,\nu}_{\theta,s\phi}
    \end{array}
    \right)
\end{displaymath}
Identify $\kernel N_{0}$ with $\C^{2}$ via the basis $\{1,\phi\}$. 
Let $N'_{s}$ denote the derivative of $N_{s}$ with respect to $s$.
Denote by
$\widehat{N}'_{0}$ the restriction of $N'_{0}$ to $\kernel N_{0}$ followed by $L^{2}$-projection
onto the same. Then from (\ref{delta-0-equ}), (\ref{delta-1-equ}), (\ref{delta-0*-equ}) we have
\begin{displaymath}
\widehat{N}'_{0}(z,w) = (-w, -z).
\end{displaymath}
This has matrix exactly as claimed in the sublemma.
This completes the proof. \qed


\subsection{Proof of Proposition~\protect\ref{struct-prop}}\label{sub-proof}

Let $W^{\alpha}_{\sigma}$, $\alpha\in\{0,I,II\}$ denote a reducible strata. Assume it is 
regular; then it is a 1-manifold which parameterizes the family of normal operators 
$N^{\alpha}_x$. (Strictly speaking we ought to work with a
lift of $W^{\alpha}_{\sigma}$ to the configuration space 
$\ZZ\times[0,1]$ but in any 
case the family will be unique up to conjugation.)

Call $W^{\alpha}_{\sigma}$ \textit{normally transverse} if given any
$x\in W^{\alpha}_{\sigma}$ there exists a 1-1 parameterization
$J\colon(-\epsilon,\epsilon)\to W^{\alpha}_{\sigma}$ of a 
neighbourhood of $x$ such that the pull-back family 
$N^{\alpha}\circ J$ is a transverse family with respect to 
spectral flow.

We need a standard result on the local structure of a bifurcation 
point which for instance  is proven in \cite{lim0}. However we 
include the proof since we shall need an ingredient from it in for
$\S\ref{sec-orient}$.

\begin{lemma}\label{bif2-lem}
    Assume $W^{\alpha}_{\sigma}$ is regular and normally transverse. 
    Let $x\in W^{\alpha}_{\sigma}$ be a singular point. Then a local 
    model for a neighbourhood $U$ of $(x,t)$ in $W_{\sigma}$ is $\Upsilon^{-1}(0)$,
    where $\Upsilon\colon\R\times[0,\infty)\to\R$, $(x,y)\mapsto xy$
    with $(0,0)\leftrightarrow x$, $\R\times\{0\}\cong U\cap 
    W^{\alpha}_{\sigma}$.
\end{lemma}

\proof $\Upsilon$ is the quadratic approximation for the Kuranishi obstruction 
map $\mXi$, below. To recall: let $\tau$ represent a tangent vector to 
$W^{\alpha}_{\sigma}$ at $(x,t)$. ($(x,t)$ is fixed throughout this 
proof.)
Then $\kernel 
\widetilde{L}_{x,t}=\R\{\tau\}\oplus\bH^{1,\nu}_{x}$ and $\coker 
\widetilde{L}_{x,t}=\bH^{1,\nu}_{x}$. Let $\Pi$ denote $L^{2}$-projection
onto $\bH^{1,\nu}_{x}$. By the implicit function theorem 
there exists a map $f\colon\R\{\tau\}\oplus\bH^{1,\nu}_{x}\to 
(\bH^{1,\nu}_{x})^{\perp}$ such that for all $r$ and $h$
\begin{eqnarray}
   &(I-\Pi)\widetilde{\XX}((t,x)+r\tau+h+f(r\tau,h))=0,\\ 
   &f(0,0)=0, \quad
   df_{0,0}=0.\nonumber
   \end{eqnarray}
   Here $\widetilde{\XX}$ denotes the maps in (\ref{para-chi}). (Note: for 
   convenience we henceforth reverse the order of the variables for 
   $\widetilde{\XX}$.)
The obstruction map 
$\mXi\colon\R\{\tau\}\oplus\bH^{1,\nu}_{x}\to\bH^{1,\nu}_{x}$ is given 
as 
\begin{equation}\label{obs-map}
\mXi(r\tau,h)=\Pi\circ\widetilde{\XX}((t,x)+r\tau+h+f(r\tau,h)).
\end{equation}
Assume, for convience that $J'(0)=\tau$.
Let $a+b\in X_{x}$ with $a\in X^{\tau}_{x}$ and $b\in X^{\nu}_{x}$. 
Then it is seen that
\begin{equation}\label{X-equ}
    \widetilde{\XX}(t',x+a+b)=K^{\alpha}_{x,t'}(a)+N^{\alpha}_{x,t'}(b)+B(a+b,a+b)
\end{equation}
where $B$ is a bilinear term. The expressions $B(a,a), B(b,b)\in 
X^{\tau}_{x}$ and $B(a,b), B(b,a)\in X^{\nu}_{x}$. The normal component of 
the second derivative of $\widetilde{\XX}$ at $(t,x)$ in the pair of directions $(r\tau,h)$
is given by
\begin{equation}
    r\frac{d}{d u}
    (N^{\alpha}\circ J(u))\bigg|_{u=0}(h)
\end{equation}
This in turn is exactly $c_{0}rh$ where $c_{0}\neq 0$ has the same sign as the spectral flow for 
$N^{\alpha}\circ J(u)$ at $u=0$ (see \cite{lim0}). 
In figuring the quadratic approximation for $\mXi$ in (\ref{obs-map}) 
we can drop the $f$ term since $df=0$ at the origin. Thus the 
quadratic approximation for $\mXi(r\tau,h)$ is $c_{0}rh$.
After dividing 
out by the action of $\stab(x)$ on $\bH^{1,\nu}_{x}$ we obtain 
$\Upsilon$. \qed

Proposition~\ref{struct-prop} largely follows directly from this 
lemma as soon as we can find an admissible perturbation $\sigma$ such that 
$W^{\alpha}_{\sigma}$ is regular and normally transverse for all 
$\alpha$, and $W^{*}_{\sigma'}$ is also regular. The remaining points
to show are the nature of the limiting ends associated to each singular 
point, and the non occurence of singular points on the $SU(3)$-trivial 
strata.

\subsubsection{Local Models}\label{sect-local-mod}

Assume $W^{\alpha}_{\sigma}$ is regular but not yet normally
transverse. Let $(x,t)$ be a point where spectral flow occurs for 
$N^{\alpha}$ and $J\colon(-\epsilon,\epsilon)\to W^{\alpha}_{\sigma}$ a
1-1 parameterization of a small neighbourhood of $(x,t)$, giving us a 
1-parameter family of operators $N_{s}=N^{\alpha}\circ 
J(s)$.

Let $\bH^{\nu}=\kernel N^{\alpha}_{x,t}$. This 
consists of $L^{2}_{3}$-sections and is independent of the completion
$N=N^{\alpha}_{x,t}\colon L^{2}_{2+k}(W^{\nu}\oplus V^{\nu})\to 
L^{2}_{1+k}(W^{\nu}\oplus V^{\nu})$, $k\ge 0$ due to the standard 
elliptic estimates. Furthermore we have $L^{2}$-Hodge-decompositions into 
closed subspaces (and using the symmetry property of $N$):
\begin{eqnarray*}
    L^{2}_{1}(W^{\nu}\oplus V^{\nu})={\rm Ran}(N|_{L^{2}_{2}})\oplus 
    \bH^{\nu}\\
    L^{2}_{2}(W^{\nu}\oplus V^{\nu})={\rm Ran}(N|_{L^{2}_{3}})\oplus 
    \bH^{\nu}.
\end{eqnarray*}
Thus $N$ by restriction,  defines a Banach space isomorphism
${\rm Ran}(N|_{L^{2}_{3}})\to {\rm Ran}(N|_{L^{2}_{2}})$. 

Assume for $J$ that the parameterization 
$N_{s}=N^{\alpha}\circ J(s)$ is at least $C^{2}$; this 
means the second derivative $\frac{d^{2}N_{s}}{ds^{2}}$ is a bounded map $L^{2}_{2}(W^{\nu}\oplus 
V^{\nu})\to L^{2}_{1}(W^{\nu}\oplus V^{\nu})$ and continuous in $s$
where the target space has the operator topology. Let $\Pi$ denote $L^{2}$-projection 
of ${\rm Ran}(N|_{L^{2}_{2}})\oplus \bH^{\nu}$ onto the first factor.
Define
\begin{displaymath}
    N^{0}_{s}= \Pi\circ N_{s}\colon{\rm Ran}(N|_{L^{2}_{3}})\to {\rm Ran}(N|_{L^{2}_{2}}).
\end{displaymath}
By the implicit function theorem the exists a subinterval, which we 
can take to be $(-\epsilon,\epsilon)$ again, and a $C^{2}$-map 
$f\colon \bH^{\nu}\times (-\epsilon,\epsilon)\to {\rm Ran}(N|_{L^{2}_{3}})$
such that the following hold:
\begin{eqnarray*}
   N^{0}_{s}(\phi+f(\phi,s)) &=&0\\
   f(\phi,0) &=& 0\\
   \frac{\partial f}{\partial s}(\phi,0)&=&0\\
   f(\cdot,s) &\in& B(\bH^{\nu},{\rm Ran}(N|_{L^{2}_{3}})).
\end{eqnarray*}
Thus solving $N_{s}(\phi)=0$, $|s|\le \epsilon$ is equivalent to 
solving the equation
\begin{displaymath}
   T(s)\phi := (I-\Pi)\circ N_{s}(\phi+f(\phi,s))=0, \quad \phi\in \bH^{\nu},
\end{displaymath}
where $T\colon(-\epsilon,\epsilon)\to{\rm Hom}(\bH^{\nu})$. It is 
clear then that $\dimn(\kernel N_{s})= \dimn(\kernel T(s))$. $T(s)$ is 
the \textit{local model} for the family $N_{s}$.

\begin{lemma}
    (a) $T(s)$ is symmetric, i.e. $\langle{T(s)\psi,\phi}\rangle_{L^{2}}
    =\langle{\psi,T(s)\phi}\rangle_{L^{2}}$ for $\psi,\phi\in 
    \bH^{\nu}$ (b) $T(s)$ commutes with $\stab(x)$.
\end{lemma}
    
Item (a) is proven by the computation
\begin{eqnarray*}
    \langle{T(s)\psi,\phi}\rangle_{L^{2}} &=& 
    \langle{N_{s}(\psi+f(s,\psi)),\phi}\rangle_{L^{2}}\\
    &=&\langle{N_{s}(\psi+f(s,\psi)),\phi+f(s,\phi)}\rangle_{L^{2}}\\
    &=& \langle{\psi+f(s,\psi),N_{s}(\phi+f(s,\phi))}\rangle_{L^{2}}\\
    &=& \langle{\psi,N_{s}(\phi+f(s,\phi))}\rangle_{L^{2}}.
\end{eqnarray*}
Item (b) follows easily from Proposition~\ref{nor-sym-prop}. Thus our 
local model is really a map into the linear space
of symmetric operators on 
$\bH^{\nu}$ which commute with $\Gamma=\stab(x)$: 
\begin{equation}\label{model-T}
    T\colon (-\epsilon,\epsilon)\to {\rm 
    Sym}_{\Gamma}(\bH^{\nu}).
\end{equation}

Let $B^{0}$ be a sufficiently 
small open neighbourhood of the origin in any finite dimensional vector 
subspace of the perturbations $\PP_{0}$ which vanish on the reducibles
(we do not want the perturbation to move $W^{\alpha}_{\sigma}$ itself).
Then we have a 
\textit{parameterized local model} 
\begin{equation}\label{model-T2}
    \widehat{T}\colon (-\epsilon/2,\epsilon/2)\times B^{0}\to {\rm 
    Sym}_{\Gamma}(\bH^{\nu})
\end{equation}
such that $\widehat{T}(s,0)=T(s)$, $s\in(-\epsilon/2,\epsilon/2)$.

\subsubsection{Completion of the argument}

\textbf{The SW-case} \textit{Step1:} start at the lowest strata
$W^{0}=\{[\mTheta]\}\times[0,1]$ which is always regular (in the 
trivial strata). We need to find an admissible perturbation in $\PP_{0}$ which 
makes $W^{0}$ normally transverse. Let $(x_{0},t_{0})\in W^{0}$ and 
let $N_{s}=N^{0}\circ J(s)$ where $J\colon(-\epsilon,\epsilon)\to W^{0}$
parameterizes a sufficiently small neighbourhood of $(x_{0},t_{0})$ 
such that a local model $s\mapsto T(s)\in\sym_{\Gamma}(\bH^{\nu})$
for $N_{s}$ is valid.

Let $\sigma\in\PP_{0}$. Then $W^{0}=W^{0}_{\sigma}$ (actually for the 
trivial strata this is valid for all $\sigma$) and $N_{s}$ changes to 
a new family $N^{\sigma}_{s}$ which is still parameterized by 
$(-\epsilon,\epsilon)$.
If $\sigma$ is sufficiently small then we 
continue to have a local model $s\mapsto \widehat{T}(s,\sigma)$, 
$s\in(-\epsilon/2,\epsilon/2)$.
In particular $N^{\sigma}_{s}$ has transverse spectral flow if and only 
if the family $s\mapsto\widehat{T}(s,\sigma)$ has transverse spectral flow. The last statement is in turn
equivalent to the following: the transformations in ${\rm 
Sym}_{\Gamma}(\bH^{\nu})$ which have non-trivial rank form a  
codimension one real subvariety, being the zeros of the determinant 
map to $\R$. Denote by $V^{(k)}$ those which 
have real
rank $\ge 4k\ge0$. Then transverse spectral flow is the condition that $s\mapsto 
\widehat{T}(s,\sigma)$ is disjoint from $V^{(k)}$, $k\ge 2$ and meets $V^{(1)}$ transversely.

If $\pi$ is an admissible perturbation on $\ZZ$ and $h$ is a function 
with support in $(0,1)$ then  $\sigma=(h\pi,0)$ is an admissible 
perturbation on $\ZZ\times [0,1]$. Let $s_{0}\in 
(-\epsilon/2,\epsilon/2)$
and $h$ have support near $s_{0}$. Then according to 
Proposition~\ref{pert-2}, there exists
$h_{i}\pi_{i}$, $i=1,\dots,n$ and $r_{i}$, $i=1,\dots,n$ sufficiently small 
such that the linearization of 
$(r_{1},\dots,r_{n})\mapsto \widehat{T}(s_{0},r_{1}h_{1}\pi_{1}+\dots +
r_{n}h_{n}\pi_{n})$ is surjective at $r_{1}=\dots=r_{n}=0$. Since 
surjectivity is an open condition this continues to hold for all $s$ 
close to $s_{0}$.
By an open cover argument we may enlarge our set
of $h_{i}\pi_{i}$ (but keeping it finite) and a sufficiently small open 
neighbourhood $B^{0}$ of the origin in their span in $\PP_{0}$ such that
$\widehat{T}\colon (-\epsilon/3,\epsilon/3)\times B^{0}\to 
\sym_{\Gamma}(\bH^{\nu})$
is submersion along $(-\epsilon/3,\epsilon/3)\times\{0\}$. Thus by 
Sard-Smale there exists an admissible perturbation $\sigma$ which makes 
$W^{0}_{\sigma}$ normally transverse along the portion
$J(-\epsilon/3,\epsilon/3)$. Another open cover argument gives normal 
transversality at all points in $W^{0}$.

\textit{Step2:} At this stage $W^{*}_{\sigma}$ and $W^{I,II}_{\sigma}$ 
are now regular in a neighbourhood of $W^{0}=W^{0}_{\sigma}$. In 
particular the local models for the bifurcation/singular points 
(Lemma~\ref{bif2-lem}) are valid for those in $W^{0}$. 
We claim that any bifurcation arc limiting into $W^{0}$ 
must come from the Type II strata. This
follows from the observation that $\bH^{1,\nu}_{\mTheta,t}=\stab(x)\cdot 
\bH^{1,II,\tau}_{\theta,t}$ at such a bifucation point. Here 
$\bH^{1,II,\tau}_{\theta,t}$ is the first 
cohomology of the normal component of the fundamental elliptic 
complex along the Type II strata (see Proposition~\ref{ellip-split}), extended over the trivial connection $\theta$.
In fact 
$\bH^{1,\nu}_{\mTheta,t}=\bH^{1,II,\tau}_{\theta,t}\oplus 
{\mathbf J}(\bH^{1,II,\tau}_{\theta,t})$ where ${\mathbf J}$ is the constant gauge 
transformation which switches the factors in the splitting 
$E=L_{0}\oplus L_{1}=(\C\oplus\C)\times Y$. 
$\bH^{1,II,\tau}_{\theta,t}\cong\R$ 
is the limiting tangent space to $W^{II}_{\sigma}$ at the end which 
limits to $([\mTheta],t)$. 

The slice space at $(x_{0},t_{0})$ for the action of $\G$ on $\ZZ\times[0,1]$ 
is $X_{x_{0}}\times\R$. Let $(v,\tau)\in X_{x_{0}}\times \R$ and $h$ a 
function with support in $(0,1)$. Then by Proposition~\ref{pert-1} there 
exists a $\pi$ such that the perturbation $(h\pi,h\tau)$ is admissible on 
$\ZZ\times[0,1]$ and $(h(t_{0})\pi_{x_{0}},h(t_{0})\tau)=(v,\tau)$. 
Furthermore this perturbation can be assumed to be supported away from $W^{0}$. 
Thus by a standard
transversality argument there exists an admissible $\sigma'$ which  
makes $W^{II}_{\sigma+\sigma'}$ regular (in 
$\CC^{II}/\G_{E}\times[0,1]$). As such it is a 1-manifold with boundary
$\M^{\sw,II}_{g_{0},\pi_{0}}\cup\M^{\sw,II}_{g_{1},\pi_{1}}$ and 
with ends (if any) which limit into $W^{0}$.

\textit{Step3:} Repeat \textit{Step1} but applied to 
$W^{II}_{\sigma+\sigma'}$. That is we can find a further perturbation 
$\sigma''\in\PP_{0}$ such that $W^{II}_{\sigma+\sigma'+\sigma''}$ is 
normally transverse. To simplify notation continue to denote by 
$\sigma$ the perturbation $\sigma+\sigma'+\sigma''$. All bifurcations 
on $W^{II}_{\sigma}$ are into $W^{*}_{\sigma}$.

\textit{Step4:} Repeat \textit{Step2} and \textit{Step3} but applied to $W^{I}_{\sigma}$.
In this \textit{Step2} since $\Gamma\cong U(1)$, $V^{(k)}$ are the 
symmetric maps which have real
rank $\ge 2k\ge0$. The rest of the argument proceeds as before.

\textit{Step5:} Repeat \textit{Step2} but applied to $W^{*}_{\sigma}$ 
so as to make it regular. This completes the proof in the SW-case.

\textbf{The SU(3)-case} We follow the same argument as above but we 
may start at \textit{Step2} as $W^{0}$ is always isolated. \qed


\section{Orientation}\label{sec-orient}

We continue to enforce the notational conventions stated at the 
beginning of \S\ref{proof-main}.

\subsection{Convention for determinant lines}\label{subsect-det}

Let $\{L_{x}\colon\mathcal{V}_{0}\to\mathcal{V}_{1}\}$ be a family of 
Fredholm operators parameterized by 
$x\in X$. The determinant line $\detind L$ is the (real) line bundle 
over $X$ whose fiber at $x$ is formally $\mLambda^{\rm max}(\kernel 
L_{x})\otimes \mLambda^{\rm max}(\coker L_{x})^{*}$. 
In the context we have been working in $L_{x}$ is a compact perturbation of a first order elliptic 
operator over $Y$ with $\mathcal{V}_{0}=L^{2}_{k+1}(V_{0})$, 
$\mathcal{V}_{1}=L^{2}_{k}(V_{1})$ and $X$ either $\CC_{E}$ or 
$\A_{F}$. 

Given a differentiable path $\gamma\colon[a,b]\to X$ we can consider
the family $L_{\gamma(t)}$ along the path $\gamma$.
An orientation of $(\detind L)_{\gamma(a)}$ can then 
be propagated to an orientation of $(\detind L)_{\gamma(b)}$ along 
$\gamma$. Since the kernel of $L_{\gamma(t)}$ may jump as $t$ varies, 
to carry out this procedure we need to stabilize the pulled back family
$\detind (\gamma^{*}L)$ over $[a,b]$. This consists of the following 
data: a finite dimensional vector space $W$ and a map $\mPsi\colon 
W\to\mathcal{V}_{1}$ such that $\widetilde{L}_{t}=L_{\gamma(t)}+\mPsi\colon 
\mathcal{V}_{0}\oplus W\to\mathcal{V}_{1}$ is surjective for all $t$.
Then $\detind (\gamma^{*}L)$ is realized as the line bundle
with fiber at $t$ being
$\mLambda^{\rm max}(\kernel \widetilde{L}_{t})\otimes
\mLambda^{\rm max}W^{*}$.

The issue we wish to address is how an orientation of a fiber of $\mLambda^{\rm max}(\kernel 
L_{\gamma(t)})\otimes \mLambda^{\rm max}(\coker L_{\gamma(t)})^{*}$ is
to be carried over to the corresponding fiber of the stabilization, 
$\mLambda^{\rm max}(\kernel \widetilde{L}_{t})\otimes
\mLambda^{\rm max}W^{*}$. Since this is a fiberwise convention we 
henceforth simplify our discussion by  considering only a single 
operator $L=L_{x}\colon\mathcal{V}_{0}\to \mathcal{V}_{1}$ in the 
family.

Suppose $\widetilde{L}=L+\mPsi\colon\V_{0}\oplus W\to \V_{1}\subset 
\V_{1}\oplus W$ is a 
stabilization. Regarding $L$ and $\widetilde{L}$ as two step chain 
complexes, we have the following exact sequence of complexes:
\begin{equation}
\begin{array}{ccccccc}
    0 & \longrightarrow & \V_{0} & \stackrel{L}{\longrightarrow} & \V_{1} & 
    \longrightarrow & 0  \\
    \downarrow &  & \downarrow &  & \downarrow &  & \downarrow  \\
    0 & \longrightarrow & \V_{0}\oplus W & 
    \stackrel{\widetilde{L}}{\longrightarrow} & \V_{1}\oplus W & \longrightarrow & 
    0  \\
    \downarrow &  & \downarrow &  & \downarrow &  & \downarrow  \\
    0 & \longrightarrow & W & \stackrel{0}{\longrightarrow} & W & \longrightarrow & 0
\end{array}
\end{equation}
This gives a long exact sequence:
\begin{equation}\label{exact-seq}
    0\to\kernel L\to \kernel\widetilde{L}\to W\to \coker L\to 
    \coker\widetilde{L}\to W\to 0.
\end{equation}
Assuming inner products on all the spaces we can \lq roll-up\rq\ the
exact sequence into a single isomorphism
\begin{equation}\label{xi-equ-2}
    \xi\colon \kernel\widetilde{L}\oplus\coker L\oplus W
    \to \kernel L\oplus W\oplus \coker\widetilde{L}.
\end{equation}
Now an orientation of $\detind L$ or $\detind \widetilde{L}$ is equivalent 
to saying that an orientation for the kernel is determined by an 
orientation of the cokernel, and vise-versa. Regarding an orientation 
$o(V)$ as a non-zero element in the highest exterior power of $V$, orient 
$\detind\widetilde{L}$ from $\detind L$ by the rule
\begin{equation}\label{xi-equ-3}
    \xi^{*}( o(\kernel L)\wedge o(W)\wedge o(\coker\widetilde{L}))
    =o(\kernel\widetilde{L})\wedge o(\coker L)\wedge o(W).
\end{equation}
The rule is independent of choice of orientatation of $W$ used as 
well as
the inner product since changes in this will only change $\xi^{*}$ by a positive 
constant.
   
We collect some remarks:
\begin{enumerate}
    \item \label{first} Allow $\mPsi$ to vary continuously with respect to 
    say a real parameter and keeping the surjectivity condition; the spaces 
    and maps in the sequence (\ref{exact-seq}) will vary continuously. 
    Thus if $\mLambda^{\max}(\kernel L)\otimes\mLambda^{\max}(\coker 
    L)^{*}$ is oriented and one stabilization is homotopic to another, 
    then the induced orientations are carried continuously onto 
    each other.
    \item \label{sec} The first statement in (\ref{first}) remains true if 
    in addition we allow $L_s$ to vary continuously with the 
    parameter $0\le s\le 1$ with the assumption the kernel does not jump and 
    $\widetilde{L}_{s}=L_{s}+\mPsi_{s}$ remains surjective.
    Therefore if $\mLambda^{\rm max}(\kernel 
    L_{\gamma(s)})\otimes \mLambda^{\rm max}(\coker L_{\gamma(s)})^{*}$ 
    is a locally trivial family then
    the propagated orientation from $s=0$ to $s=1$ is consistent with 
    the propagated orientation after stabilzation.
    \item It is not neccesary to insist on stabilizations which are 
    \lq surjective\rq\ in (\ref{first}) and (\ref{sec}). 
    For instance everything we 
    have said so far is equally true applied to say the stabilization 
    where $\mPsi$ is the zero map. The homotopy 
    assertions  continue to be true with the assumption the kernel does not
    jump in a continuous variation.
\end{enumerate}

\subsection{The Floer-Taubes operator in the parameterized 
context}\label{sect-FT}

This subsection contains the proof of Proposition~\ref{orient-para-lem}.

Retain the notation of \S\ref{proof-main}. 

\lq Rolling up\rq\ the complex (\ref{para-ellip}) gives a (perturbed) elliptic operator
\begin{equation}
    \widetilde{L}_{x,t}\colon L^{2}_{2}(W\oplus V)\oplus\R
    \to L^{2}_{1}(W\oplus V).
\end{equation}
This operator has the property that the restriction to $L^{2}_{2}(W\oplus 
V)$ for a fixed value of $t$ gives the Floer-Taubes operator on $Y$
with respect to parameter value $t$. As in \S\ref{splitting}, along a 
reducible strata $\ZZ^{\alpha}\times[0,1]$ we can take the
 \textit{tangential} component of $\widetilde{L}_{x,t}$:
\begin{equation}
    \widetilde{K}^{\alpha}_{x,t}\colon L^{2}_{2}(W^{\tau}\oplus V^{\tau})\oplus\R
    \to L^{2}_{1}(W^{\tau}\oplus V^{\tau}).
\end{equation}
The determinant index $\detind \widetilde{L}_{x,t}$ 
determines an orientation of $W^{*}_{\sigma}$ and $\detind 
\widetilde{K}^{\alpha}$ and orientation of $W^{\alpha}_{\sigma}$. 
An orientation of $\detind \widetilde{L}_{x,t}$ and $\detind \widetilde{K}^{\alpha}$
will immediately give Proposition~\ref{orient-para-lem} up to an overall 
sign in the various stratas. The issue is how to fix the overall orientation by  
so that the sign is correct. This in turn reduces to fixing the orientation along the trivial strata.

For convenience, denote either $\widetilde{L}_{x,t}$ or 
$\widetilde{K}^{\alpha}_{x,t}$ by $\widetilde{\mLambda}_{x,t}$, and 
$\widetilde{\mLambda}^{0}_{x,t}$ the restriction to the summand different 
from $\R$. At a point on the 
trivial strata, say represented by $(\mTheta,t)$, the kernel and 
cokernel ($=$ $L^{2}$-orthogonal of range) of $\widetilde{\mLambda}^{0}$ coincide by self-adjointness.
Thus we have
\begin{equation}
    \kernel \widetilde{\mLambda}_{\mTheta,t}=\kernel 
    \widetilde{\mLambda}^{0}_{\mTheta,t}\oplus\R, \quad
    \coker\widetilde{\mLambda}_{\mTheta,t}=\coker\widetilde{\mLambda}^{0}_{\mTheta,t}
    =\kernel \widetilde{\mLambda}^{0}_{\mTheta,t}.
\end{equation}
The $\R$ summand in the kernel corresponds exactly to the tangent 
space to the trivial strata $\{[\mTheta]\}\times[0,1]$.

\begin{prop}
    Orient $\widetilde{\mLambda}$ by specifing the orientation 
    of $\detind \widetilde{\mLambda}$ at  $(\mTheta,t)$ by
    $\overline{o}(\R)\wedge o(\kernel \widetilde{\mLambda}^{0}_{\mTheta,t})
	\wedge o({\kernel \widetilde{\mLambda}^{0}_{\mTheta,t}})^{*}$ where
	$\overline{o}(\R)$ is the standard orientation.
	Then Proposition~\ref{orient-para-lem} holds with this choice.
\end{prop}

\proof It suffices to consider irreducible strata of the constant family 
$W_{\sigma}=\M_{\pi}\times[0,1]$; the reducible stratas are treated in 
an identical manner using the tangential operators for that strata.

Let $[x]\in\M_{\pi}^{*}$ and $\gamma\colon[0,1]\to\ZZ\times\{t\}$ a path from 
$(\mTheta,t)$ to $x$. The assumption that $W_{\sigma}$ is a constant 
family gives 
\begin{equation}\label{kernel-equ}
    \kernel \widetilde{\mLambda}_{\gamma(t)}=\kernel 
    \widetilde{\mLambda}^{0}_{\gamma(t)}\oplus\R, \quad
    \coker\widetilde{\mLambda}_{\gamma(t)}
    =\kernel \widetilde{\mLambda}^{0}_{\gamma(t)}.
\end{equation}
Note that $\widetilde{\mLambda}^{0}$ is the family which orients 
$\M_{\pi}^{*}$. Equation (\ref{kernel-equ}) gives us an isomorphism 
from $\detind (\gamma^{*}\widetilde{\mLambda}^{0})$ to 
$\detind (\gamma^{*}\widetilde{\mLambda})$ by the rule
$\omega\mapsto \overline{o}(\R)\wedge\omega$. This isomorphism 
carries the standard orientation for $\detind\widetilde{\mLambda}^{0}$
at $\mTheta$ to the orientation of $\detind\widetilde{\mLambda}^{0}$
at $(\mTheta,t)$ given in the statement of the proposition. Thus if 
$\epsilon(x)$ is the orientation of $x$ in $\ZZ^{*}/\G$ then the arc
$\{x\}\times [0,1]=[0,1]$ in $W^{*}_{\sigma}$ is assigned the 
orientation $\epsilon(x)\overline{o}(\R)$. Now it is straightforward
to see that $\partial 
W^{*}_{\sigma}=\M^{*}_{\pi}\times\{1\}-\M^{*}_{\pi}\times\{0\}$. \qed

\subsection{Bifurcation Points}

This subsection estabishes Proposition~\ref{sing-pt-lem}.

\subsubsection{Type I, II Strata}

Let us assume that the various strata of $W_{\sigma}$ are oriented 
according to the convention in the preceding section. Let $(x,t_{0})$ 
denote a bifurcation point on a Type I or II strata; parameterize a 
small neighbourhood by $J\colon(-\epsilon,\epsilon)\to W_{\sigma}$ 
with $J(0)=(x,t_{0})$ and consistent with the orientation of the 
reducible strata.
Then the family of normal 
operators $N_{u}=N^{\alpha}\circ J(u)$ has transverse spectral flow at $u=0$ and $\kernel 
N_{0}\cong \C\cong\coker N_{0}$.
We have $\stab(x)\cong U(1)$ so $\bH^{0}_{x}\cong\R$, and $\kernel 
\widetilde{L}^{\tau}_{x}=\R\{\tau\}\oplus\bH^{0}_{x}$
where $\tau\in\widetilde{\bH}^{1,\tau}_{x}$ represents a non-zero tangent vector to the 
strata at $x$. We assume that $\tau$ is consistent with the 
orientation of the strata at $x$.

Let $v\in\bH^{1}_{x}$ be unit length. To 1st order, the bifurcation 
arc at $x$ is modelled by $\{sv\}$, $0\le s<\epsilon$. Now define a 
1-parameter family $\widetilde{L}_{s}=\widetilde{L}_{\xi(s)}$ where 
$\xi$ parametrizes (up to gauge equivalence) the bifurcation arc with $\xi(0)=x$ and $\xi'(0)=v$. 
We have
\begin{equation}
     \kernel\widetilde{L}_{0}=\R\{\tau\}\oplus\bH^{0}_{x}\oplus
     \bH^{1}_{x},
     \quad
     \coker\widetilde{L}_{0}=\bH^{0}_{x}\oplus\bH^{1}_{x}
\end{equation}
where we identify the cokernel with the $L^{2}$-orthogonal of the 
range. On the other hand for $s\neq 0$, $\kernel \widetilde{L}_{s}\cong \R$ is 
the tangent space to the bifurcation arc, whilst $\coker \widetilde{L}_{s}=\{0\}$.
Notice that $\lim_{s\to 0}\kernel \widetilde{L}_{s}=\R\{ v\}$.

Let $\gamma\in\bH^{0}_{x}$ be unit length. This is a tangent 
vector to $\stab(x)$ at the identity.
The derivative of the $\stab(x)$ action on $\bH^{1}_{x}$ gives an 
action of $\gamma$ on $\bH^{1}_{x}$ which is a complex structure. We 
fix the complex structure by our choice of $\gamma$. Using $\tau$, 
$\gamma$ and $v$ as a basis we may express
\begin{equation}\label{basis-equ}
   \kernel\widetilde{L}_{0}=\R\{\tau\}\oplus\R\{\gamma\}
   \oplus\C\{v\},
   \quad 
   \coker\widetilde{L}_{0}=\R\{\gamma\}
   \oplus\C\{v\}.
   \end{equation}
   
Our next goal will be to reduce $\widetilde{L}_{s}$ to a more 
managable finite dimensional family $T_{s}$. The following discussion  
parallels subsection~\ref{sect-local-mod} so we shall be 
brief.
$\widetilde{L}_{s}$ is a bounded operator $L^{2}_{2+k}(W\oplus V)\oplus\R\to 
L^{2}_{1+k}(W\oplus V)$, $k\ge 0$. The $L^{2}$-adjoint is a bounded 
operator $\widetilde{L}^{*}_{s}\colon L^{2}_{2+k}(W\oplus V)\to 
L^{2}_{1+k}(W\oplus V)\oplus\R$. We have the following 
$L^{2}$-decompositions:
\begin{eqnarray*}
    L^{2}_{1}(W\oplus V) &=& \ran (\widetilde{L}_{0}|_{L^{2}_{2}})\oplus 
    \kernel(\widetilde{L}^{*}_{0}|_{L^{2}_{1}})\\
    L^{2}_{2}(W\oplus V) &=& \ran (\widetilde{L}^{*}_{0}|_{L^{2}_{3}})\oplus 
    \kernel(\widetilde{L}_{0}|_{L^{2}_{2}}).
\end{eqnarray*}
Let ${\cal V}=\kernel (\widetilde{L}_{0}|_{L^{2}_{3}})$, ${\cal 
W}=\kernel(\widetilde{L}^{*}_{0}|_{L^{2}_{2}})$. Denote by $\Pi_{\cal V}$, 
$\Pi_{\cal W}$ $L^{2}$-projection maps onto ${\cal V}$, ${\cal W}$ 
respectively. By the implicit function theorem the following holds:
\begin{enumerate}
    \item for $|s|<\epsilon$ there exists $f\colon{\cal 
    V}\times(-\epsilon,\epsilon)\to {\cal V}^{\perp}$ such that
    \begin{eqnarray*}
	(I-\Pi_{\cal W})\circ \widetilde{L}_{s}(v+f(v,s))&=&0\\
	f(v,0)&=& 0\\
	\frac{\partial f}{\partial s}(v,0)&=&0\\
	f(\cdot,s)&\in &B({\cal V}, \ran(\widetilde{L}^{*}_{0}|_{L^{2}_{3}}))
	\end{eqnarray*}
	\item for $|s|<\epsilon$ there exists $g\colon{\cal W}\times 
	(-\epsilon,\epsilon)\to {\cal W}^{\perp}$ such that
	\begin{eqnarray*}
	(I-\Pi_{\cal V})\circ \widetilde{L}^{*}_{s}(w+g(w,s))&=&0\\
	g(w,0)&=& 0\\
	\frac{\partial g}{\partial s}(w,0)&=&0\\
	g(\cdot,s)&\in &B({\cal W}, \ran(\widetilde{L}_{0}|_{L^{2}_{2}}))
	\end{eqnarray*}
\end{enumerate}
These lead to the following isomorphisms:
\begin{eqnarray}
    &&\gamma_{1,s}\colon L^{2}_{2}(W\oplus V)\to L^{2}_{2}(W\oplus V)\\
    &&(v,v') \mapsto (v,v'+f(v',s)), \quad v\in{\cal V}, v'\in{\cal 
    V}^{\perp}\nonumber
    \end{eqnarray}
\begin{eqnarray}
    &&\gamma_{2,s}\colon L^{2}_{1}(W\oplus V)\to L^{2}_{1}(W\oplus V)\\
    &&(w,w') \mapsto (w,w'+f(w',s)), \quad w\in{\cal W}, w'\in{\cal 
    W}^{\perp}\nonumber
    \end{eqnarray}

The preceding proves:

\begin{lemma}
    The family
    $\widetilde{L}^{0}_{s}=\gamma^{-1}_{2,s}\circ \widetilde{L}_{s}\circ 
    \gamma_{1,s}$
    is isomorphic to $\widetilde{L}_{s}$. Define
$T\colon(-\epsilon,\epsilon)\to {\rm Hom}({\cal V},{\cal W})$ by the 
rule 
    \begin{equation}\label{TL-equ}
    T_{s}(v)= \Pi_{\cal W}\circ \widetilde{L}_{s}(v+f(v,s)).
\end{equation}
    Then
    $\widetilde{L}^{0}_{s}$ splits orthogonally as $T_{s}\oplus 
    \widetilde{L}^{1}_{s}\colon{\cal V}\oplus {\cal V}^{\perp}\to{\cal 
    W}\oplus{\cal W}^{\perp}$ where $\widetilde{L}^{1}_{s}$ is an isomorphism.
    Lastly $\detind \widetilde{L}\equiv\detind \widetilde{L}^{0}\equiv\detind T$.
\end{lemma}

The last statement of the lemma is a consequence of the observation
$\kernel \widetilde{L}_{s}=\kernel\widetilde{L}^{0}_{s}=\kernel T_{s}$
and $\coker \widetilde{L}_{s}=\coker\widetilde{L}^{0}_{s}=\coker T_{s}$.
$T_{s}$ is the \textit{local model} for $\widetilde{L}_{s}$. Clearly $\detind T$ 
inherits the orientation of $\detind \widetilde{L}$.

With the aim of determining the orientation of the bifurcation arc we 
may now exclusively work with $T_{s}$ rather then $\widetilde{L}_{s}$.
For $s\neq 0$, 
$\kernel T_{s}\cong\R$ and $\coker T_{s}=\{0\}$ with $\lim_{s\to 0}
\kernel T_{s}=\R\{ v\}$. An orientation of $\lim_{s\to 0}
\kernel T_{s}$ mutually determines an orientation of $\kernel T_{s}$ for $s\neq 
0$ by continuity. Recall that $\lim_{s\to 0}\kernel T_{s}=\R\{v\}$ models (to first 
order) the bifurcation arc at $([x],t_{0})$. 
We wish to get the induced orientation on $\lim_{s\to 0}\kernel T_{s}=\R\{v\}$ 
from $\detind T$. Without loss of generality we can replace $T_{s}$ 
with the linearized family
\begin{equation}
    \widehat{T}_{s}=T_{0}+sT'_{0} =sT'_{0}
\end{equation}
since $T_{0}=0$.
$\detind\widehat{T}$ agrees with $\detind T$ at $s=0$ so it is 
oriented according to our coventions (\S\ref{sect-FT}).

\begin{lemma}\label{bif-linear}
   Identify $\kernel T_{0}=\kernel\widetilde{L}_{0}$ with $\R\oplus\R\oplus\C$
   and $\coker T_{0}=\coker \widetilde{L}_{0}$ with $\R\oplus\C$ using the basis 
   in (\ref{basis-equ}). Let $T'_{0}$ denote the derivative of 
    $T_{s}$ at $s=0$. 
    Then
    \begin{equation}\label{T2-equ}
	T'_{0}(t,r,z)=( {\rm Re}(i\overline{z}),c_{0}t+ir)
	\end{equation}
	where $c_{0}\neq 0$ has the same sign as the spectral flow for
	$N_{u}$ at $u=0$.
\end{lemma}

This will be proven below.

The proof of Proposition~\ref{sing-pt-lem} in the Type I,II cases is 
now a direct consequence of 
the next lemma. 
For the linearized family, $\kernel \widehat{T}_{s}=\R\{v\}$
 and $\coker \widehat{T}_{s}=\{0\}$, $s\neq 0$.

\begin{lemma}
Give $\R\{v\}=\R$ the standard orientation 
$\overline{o}(\R)$ of $\R$. (With this 
orientation the bifurcation arc is pointing away from the reducible 
strata.) Let $o'(\R)$ be the induced orientation on 
$\R\{v\}=\R$ given by the orientation of $\detind \widehat{T}$. Then 
$o'(\R)=-{\rm sign}(c_{0})\overline{o}(\R)$.
\end{lemma}

\proof We shall stabilize the family 
$\widehat{T}$ explicitly and evaluate the propagated orientation on 
$\detind\widehat{T}$, $s\neq 0$. Let $\V=\R\oplus\R\oplus\C$ and 
$\W=\R\oplus\C$. Stabilize $\widehat{T}_{s}$ by 
\begin{equation}
    \widehat{T}_{s}+I\colon \V\oplus \W\to\W\oplus\W.
\end{equation}
Clearly this is surjective for all $s$. Let $\V_{s}=\kernel (\widehat{T}_{s}+I)$.
Then $\V_{0}=\V$ and $\{\V_{s}\}$ is a locally trivial family. Denote 
by $p_{1}$ the projection onto the 1st factor of $\V\oplus\W$. Then 
$p_{1}\colon \V_{s}\to \V$ is an isomorphism. 
Choose an orientation $o(\W)$ for $\W$ (the result we obtain will be 
seen to be independent of the choice). 
The reference orientation of $\detind\widehat{T}$ at $s=0$ dictates 
that the induced orientation $o(\V)$ is $\tau\wedge o(\W)$.
The orientation propagated to $\V_{s}$ is obtained by 
pulling back via $p_{1}$; denote this as $o(\V_{s})=p_{1}^{*}o(\V)$. Thus 
$\detind (\widehat{T}+I)$ at $s$ has the propagated orientation 
$p_{1}^{*}o(\V)\wedge o(\W^{*})$.

On the other hand, since $\widehat{T}+I$ is a stabilization there is, 
according to our convention of \S\ref{subsect-det} a canonical way of 
relating the orientation on $\R\{v\}=(\detind \widehat{T})_{s}$ with
$\mLambda^{\rm max}\V_{s}\otimes\mLambda^{\rm max}\W^{*}=(\detind(\widehat{T}+I))_{s}$, $s\neq 0$.
This is given by the exact sequence (\ref{exact-seq}) and the rule 
(\ref{xi-equ-3}). In our situation since $\coker 
\widehat{T}_{s}=\{0\}$, $s\neq 0$, we have for $s\neq 0$, the exact 
sequence
\begin{equation}
    0\longrightarrow
    \R\{v\}\stackrel{\iota}{\longrightarrow}
    \V_{s}\stackrel{p_{2}}{\longrightarrow}
    \W\longrightarrow
    0\longrightarrow
    \W\stackrel{\rm Id}{\longrightarrow}
   \W\longrightarrow 0.
\end{equation}
Here $\iota$ is the inclusion and $p_{2}$ is projection onto the 2nd 
factor of $\V\oplus\W$. Once we have chosen $o(\W)$ as above, our 
stabilization convention requires that $\V_{s}$ has the induced 
orientation 
\begin{equation}\label{orient-equ}
o'(\V_{s})=v\wedge p_{2}^{*}o(\W)=\epsilon o(\V_{s}),
\end{equation}
$\epsilon\neq 0$. The lemma is proven once we can show 
${\rm sign}(\epsilon)=-{\rm sign}(c_{0})$. 

The composition $p_{2}\circ (p_{1})^{-1}\colon \V\to\W$ is given by 
$u\mapsto(u,-sT'_{0}(u))\mapsto -sT'_{0}(u)$. If we identify $\V_{s}$ with $\V$ 
via $p_{1}$ then (\ref{orient-equ}) is equivalent to 
$o'(\V)=-v\wedge (sT'_{0})^{*}o(\W)=\epsilon o(\V)$. Choose say $o(\W)=\gamma\wedge v\wedge iv$.
Then according to Lemma~\ref{bif-linear}, $(T'_{0})^{*}(\gamma\wedge v\wedge 
iv)=\frac{1}{c_{0}}iv\wedge\tau\wedge\gamma$. Thus 
\begin{equation}
o'(\V)=-\frac{1}{s^{3}c_{0}}v\wedge iv\wedge\tau\wedge\gamma=
-\frac{1}{s^{3}c_{0}}\tau\wedge o(\W)=-\frac{1}{s^{3}c_{0}}o(\V).
\end{equation}
Hence $\epsilon=-{\rm sign}(c_{0})$.\qed

\proofof{Lemma~\ref{bif-linear}}
It suffices to replace $\widetilde{L}_{s}$ with the \lq 1st order 
bifurcation family\rq\ 
$\widehat{L}_{s}=\widetilde{L}_{x+sv}$ since the computation takes 
place at $s=0$. In the rule (\ref{TL-equ}) with $\widetilde{L}$ 
replaced by $\widehat{L}$, we may drop the $f$ term when computing the
derivative of the family at $s=0$ since $\frac{\partial f}{\partial 
s}|_{s=0}=0$.
Then restricted to $\R\{\tau\}\oplus 
\R\{\gamma\}\oplus\C\{v\}$ we have
\begin{equation}
    \widehat{L}_{s}(t\tau,r\gamma,zv)= (\delta^{0*}_{x+sv}(zv),
    \widetilde{\delta}^{1}_{x+sv}(t\tau,zv)+\delta^{0}_{x+sv}(r\gamma)).
\end{equation}
The derivative at $s=0$ of $\delta^{0}_{x+sv}(r\gamma)$ is 
$r\gamma\cdot v=irv$ which by our convention $\gamma\cdot v$ defines the complex 
structure $i$ on $\C\{v\}$. The derivative at $s=0$ for 
$\delta^{0*}_{x+sv}(zv)$ is then the adjoint of this complex 
multiplication and hence the real part of 
$\langle{i,z}\rangle\gamma=i\overline{z}\gamma$. The last term in 
(\ref{T2-equ}) unaccounted for is $c_{0}tv$ (recall $v$ is taken as 
one of the basis vectors) which is the derivative at $s=0$ for 
$\widetilde{\delta}^{1}_{x+sv}(t\tau,zv)$. This computation essentially 
reduces to that of the quadratic approximation for the obstruction map
in the proof of Lemma~\ref{bif2-lem}. Following the notation introduced
there and equation (\ref{X-equ}) we obtain
\begin{eqnarray*}
    \widetilde{\delta}^{1}_{t',x+a+b}(t\tau,a'+b')
    &=&t\frac{d}{d u}(K^{\alpha}\circ J(u))\bigg|_{u=0}(a)
    +K_{t',x}(a')\\ 
    &&\mbox{}+t\frac{d}{d u}(N^{\alpha}\circ J(u))\bigg|_{u=0}(b)
    +N^{\alpha}_{t',x}(b') \\
    &&\mbox{}+B(a+b,a'+b')+B(a'+b',a+b).
\end{eqnarray*}
Now replacing $b$ with $sv$ and taking the derivative with respect to $s$ 
at $s=0$ gives
\begin{equation}
    \frac{d}{d s}(\Pi_{\W}\circ\widetilde{\delta}^{1}_{t',x+sb})\bigg|_{s=0}(t\tau,b')
    =t\frac{d}{d u}(N^{\alpha}\circ J(u))\bigg|_{u=0}(v)=c_{0}tv.
\end{equation}
This completes the computation of $T'_{0}$ noting that in 
(\ref{TL-equ}) we can drop the $f$ term since $df=0$ at 
$(0,0)$. \qed

\subsection{Trivial strata}

Recall that in the SW-case, any bifurcation along the trivial strata 
ends up in the Type II strata. Fix a splitting $E=L_{0}\oplus L_{1}$ 
and let $\CC(L_{0})\subset\CC_{E}$ denote the subset consisting of 
pairs $(A,\mPhi)$ where $A=A_{0}\oplus\theta$ ($\theta$ is a fixed 
trivial connection on $L_{1}$) and $\mPhi=(\phi_{0},0)$
with respect to the splitting. It is clear that any Type II reducible 
is gauge equivalent to one in $\CC(L_{0})$.

Let $\G^{0}$ be the subgroup of gauge 
transformations which take block diagonal form
$$
	\left(
	\begin{array}{cc}
	g&0\\
	0& 1
	\end{array}
	\right).
$$
$\G^{0}$ preserves 
$\CC(L_{0})$ setwise (but is not the largest subgroup to do so). 
$\CC(L_{0})$ with the action of $\G^{0}$ is basically the context of 
$U(1)$-SW theory and we
may regard the bifurcation we are concerned with as happening along the trivial strata
$\{[\theta]\}\times[0,1]$ inside 
$\CC(L_{0})/\G^{0}\times[0,1]$.

At any $(A_{0},\phi_{0},\theta)\in\CC(L_{0})$ the tangential operator $K^{\sw,II}$ splits as 
$K^{0}_{A_{0},\phi_{0}}\oplus K$ where $K$ is our (untwisted) deRham operator on 
$\mLambda^{0+1}$. Thus we see that $\detind K^{\sw,II}=\detind K^{0}$.
We may regard $K^{0}$ as being the SW-Floer-Taubes operator for 
$\CC(L_{0})$ and therefore the 
deformation theory, orientations, spectral flow, etc. are completely 
determined in reference to $K^{0}$ which in turn is determined by 
$K^{\sw,II}$. In our $U(1)$ gauge theory perspective the trivial strata has 
stabilizer $U(1)\subset\G^{0}$ and therefore we are back in the situation
of the Type I,II bifurcations. In particular the proof of the 
orientation of the bifurcation arc now extends to this situation.

\normalsize

\noindent{\it Department of Mathematics\newline
University of California\newline
Santa Cruz, CA 95064}\newline
\noindent{\tt ylim@math.ucsc.edu}

\end{document}